\UseRawInputEncoding
\documentclass[11pt]{article}
\usepackage{amsmath,amssymb,amsthm,amsfonts}
\usepackage{graphicx}
\usepackage{float}
\usepackage[scriptsize,nooneline]{subfigure}
\usepackage{indentfirst}
\usepackage{cite}
\usepackage{color}
\usepackage{mathrsfs}
\usepackage{epstopdf}
\usepackage{setspace}
\usepackage[labelfont=footnotesize,font=footnotesize]{caption}
\usepackage[colorlinks, citecolor=red]{hyperref}

\newtheorem{theorem}{Theorem}[section]

\newtheorem{lemma}{Lemma}[section]

\newtheorem{definition}{Definition}[section]

\usepackage{hyperref}
\def\geq{\geqslant}\def\leq{\leqslant}

\begin{document}
\title{\bf Characterizations of the mixed central Campanato spaces via the commutator operators of Hardy type}
\author{Wenna Lu and Jiang Zhou\thanks{The research is supported by the NNSF of China (No. 12061069).} \\[.5cm]}
\date{}
\maketitle
{\bf Abstract.} The purpose of this paper is to establish some characterizations of mixed central Campanato space $\mathfrak{C}^{\vec{p},\lambda}(\mathbb{R}^{n})$, via the boundedness of the commutator operators of Hardy type. Unlike the case $\lambda\geq0$, there are some technical difficulties caused by $\lambda<0$ to be overcome. In addition, an extra assumption called as the mixed version of the reverse
H\"{o}lder class be required in the proof of the converse characterization. Moreover, some further interesting conclusions for the Hardy type operators on mixed $\lambda$-central Morrey space $\mathcal{B}^{\vec{p},\lambda}(\mathbb{R}^{n})$ are also derived.


{\bf MSC(2020) subject classification:} 42B20, 42B25, 42B35.

\baselineskip 15pt
\section{Introduction and main results}
\quad Let $f$ be a non-negative integrable function on $\mathbb{R}^{+}$. The classical Hardy operator and its dual operator are defined by
$$Hf(x):=\frac{1}{x}\int^{x}_{0}f(t)dt, \qquad H^{*}f(x):=\int^{\infty}_{x}\frac{f(t)}{t}dt, \quad x>0.$$
Then Hardy operators $H$ and $H^{*}$ are adjoint mutually, that is,
$$\int^{\infty}_{0}Hf(x)g(x)dx=\int^{\infty}_{0}f(x)H^{*}g(x)dx,$$
when $f\in L^{p}(\mathbb{R}^{+}), g\in L^{q}(\mathbb{R}^{+}), 1<p<\infty, 1/p+1/q=1$.

In 1920, Hardy \cite{H1} first introduced the classical Hardy operator, and established the celebrated Hardy's integral inequalities as follows:
$$\|Hf\|_{L^{p}(\mathbb{R}^{+})}\leq\frac{p}{p-1}\|f\|_{L^{p}(\mathbb{R}^{+})}, \qquad \|H^{*}f\|_{L^{q}(\mathbb{R}^{+})}\leq\frac{p}{p-1}\|f\|_{L^{q}(\mathbb{R}^{+})}.$$
Moreover,
$$\|H\|_{L^{p}(\mathbb{R}^{+})\rightarrow L^{p}(\mathbb{R}^{+})}=\|H^{*}\|_{L^{q}(\mathbb{R}^{+})\rightarrow L^{q}(\mathbb{R}^{+})}=\frac{p}{p-1},$$
where $1/p+1/q=1$.

In 1976, the classical Hardy operator was extended to $n$-dimensional form by Faris \cite{FWG}, that is, let $f$ be a locally integrable function in $\mathbb{R}^{n}$, the $n$-dimensional Hardy operator $\mathcal{H}$ is defined by
$$\mathcal{H}f(x)=\frac{1}{|x|^{n}}\int_{|t|<|x|}f(t)dt, \quad x\in\mathbb{R}^{n}\backslash\{0\}.$$
The norm of $\mathcal{H}$ on $L^{p}(\mathbb{R}^{n})$ was evaluated in \cite{CG} and found to be equal to that of the one-dimensional Hardy operator in \cite{HLP}.
Similarly, the dual operator $\mathcal{H}^{*}$ of $\mathcal{H}$ is given as
$$\mathcal{H}^{*}f(x)=\int_{|t|\geq|x|}\frac{f(t)}{|t|^{n}}dt, \quad x\in\mathbb{R}^{n}\backslash\{0\}.$$
Also, it is not hard to see that $\mathcal{H}$ and $\mathcal{H}^{*}$ satisfy
$$\int_{\mathbb{R}^{n}}\mathcal{H}f(x)g(x)dx=\int_{\mathbb{R}^{n}}f(x)\mathcal{H}^{*}g(x)dx,$$
when $f\in L^{p}(\mathbb{R}^{n}), g\in L^{q}(\mathbb{R}^{n}), 1<p<\infty, 1/p+1/q=1$.

Let $b$ be a locally integrable function on $\mathbb{R}^{n}$ and  $T$ be a Calder\'{o}n-Zygmund integral operator. Then the commutator operator defined for a proper function $f$ can be denoted by
$$T_{b}(f)=:bTf-T(bf).$$
The function $b$ is called the symbol function of $T_{b}$. The investigation of the operator $T_{b}$ begins with Coifman-Rochberg-Weiss \cite{CRW} pioneering study of the operator $T$. A well-known result of Coifman et al. stated that the commutator $T_{b}f$ is bounded on $L^{p}(\mathbb{R}^{n})$ if and only if the symbol function $b$ is in the bounded mean oscillation space $\mathrm{BMO}(\mathbb{R}^{n})$.

The commutators of $\mathcal{H}$ and $\mathcal{H}^{*}$ are defined by
$$\mathcal{H}_{b}f=b\mathcal{H}f-\mathcal{H}(fb),\quad \mathcal{H}^{*}_{b}f=b\mathcal{H}^{*}f-\mathcal{H}^{*}(fb).$$
An important work on the commutator of Hardy operator  $H_{b}$  is due to Long and Wang \cite{LW}, who showed that the commutators of 1-dimensional Hardy operator $H_{b}$ and its dual operator $H^{*}_{b}$ are bounded on Lebesgue spaces if and only if the symbol function $b$ belong to a one sided dyadic central bounded mean oscillation space $\mathrm{CMO}^{p}(\mathbb{R}^{+})$. In 2007, Fu et al. \cite{FLW} gave the characterization of the central bounded mean oscillation space $\mathrm{CMO}^{p}(\mathbb{R}^{n})$ via the boundedness of the commutator of $n$-dimensional Hardy type operators $\mathcal{H}_{b}$ and $\mathcal{H}^{*}_{b}$ on Lebesgue spaces. These known facts indicate that the commutators of Hardy type operators are fundamentally different from other some important operators in harmonic analysis. The difference happens because the Hardy operator is an average operator around the origin. Therefore, this topic on the research of Hardy type operators and their commutators is very necessary and meaningful.

 Assume that $1\leq p<\infty$ and $-\frac{1}{p}<\lambda<\frac{1}{n}$. Then the central Campanato space $\mathcal{C}^{p,\lambda}(\mathbb{R}^{n})$ \cite{SL3} is defined by
$$\mathcal{C}^{p,\lambda}(\mathbb{R}^{n}):=\{f\in L^{p}_{loc}(\mathbb{R}^{n}):\|f\|_{\mathcal{C}^{p,\lambda}(\mathbb{R}^{n})}<\infty\},$$
where
$$\|f\|_{\mathcal{C}^{p,\lambda}(\mathbb{R}^{n})}:=
\sup\limits_{r>0}\frac{1}{|B(0,r)|^{\lambda}}\Big(\frac{1}{|B(0,r)|}\int_{B(0,r)}|f(x)-f_{B(0,r)}|^{p}dx\Big)^{\frac{1}{p}},$$
here and in the sequel,
$$f_{B(0,r)}:=\frac{1}{|B(0,r)|}\int_{B(0,r)}f(t)dt.$$
If $\lambda=0$, it is not hard to see that $\mathcal{C}^{p,0}(\mathbb{R}^{n})=CMO^{p}(\mathbb{R}^{n})$, this space was introduced by Lu and Yang \cite{LY} with the equivalent norm
$$\|f\|_{\mathrm{CMO}^{p}(\mathbb{R}^{n})}:=
\sup\limits_{r>0}\inf\limits_{c\in\mathbb{C}}\Big(\frac{1}{|B(0,r)|}\int_{B(0,r)}|f(x)-c|^{p}dx\Big)^{\frac{1}{p}}.$$
In fact, the space $CMO^{p}(\mathbb{R}^{n})$ can be regarded as a local version of $BMO(\mathbb{R}^{n})$ at the origin. However, their properties are quite different; for example, the famous John-Nirenberg inequality for the space $BMO(\mathbb{R}^{n})$ no longer holds in the space $CMO^{p}(\mathbb{R}^{n})$. In addition, if $0<\lambda<\frac{1}{n}$, then $\mathcal{C}^{p,\lambda}(\mathbb{R}^{n})$ is the $\lambda$-central bounded mean oscillation space $CMO^{p,\lambda}(\mathbb{R}^{n})$, which was introduced by Alvarez,  Guzm\'{a}n-Partida and Lakey in \cite{AL}. Furthermore, we have $\mathcal{C}^{p,\lambda}(\mathbb{R}^{n})\supset\mathcal{B}^{p,\lambda}(\mathbb{R}^{n})$ when $-\frac{1}{p}<\lambda<0$. Here, $\mathcal{B}^{p,\lambda}(\mathbb{R}^n)$ denotes the central Morrey spaces \cite{AL} with the following norm
$$\mathcal{B}^{p,\lambda}(\mathbb{R}^n):=\{f\in L^{p}_{loc}(\mathbb{R}^{n}):\|f\|_{\mathcal{B}^{p,\lambda}(\mathbb{R}^n)}<\infty\},$$
where
$$\|f\|_{\mathcal{B}^{p,\lambda}(\mathbb{R}^{n})}:=
\sup\limits_{r>0}\frac{1}{|B(0,r)|^{\lambda}}\Big(\frac{1}{|B(0,r)|}\int_{B(0,r)}|f(x)|^{p}dx\Big)^{\frac{1}{p}}.$$

For the case $0<\lambda<\frac{1}{n}$, Zhao and Lu \cite{ZL} proved the following result.
\begin{theorem}[\cite{ZL}]
Let $1<p, q<\infty$, $0\leq\lambda<\frac{1}{n}$ and $\lambda=\frac{1}{p}-\frac{1}{q}$. Then both $\mathcal{H}_{b}$ and $\mathcal{H}^{*}_{b}$ are bounded from $L^{p}(\mathbb{R}^n)$ to $L^{q}(\mathbb{R}^n)$ if and only if
$$b\in\mathrm{CMO}^{\max\{q,p'\},\lambda}(\mathbb{R}^{n}).$$
\end{theorem}

For the case $-\frac{1}{p}<\lambda<0$, Shi and Lu \cite{SL3} obtained the corresponding conclusions under an additional assumption, the well-known mean value inequality, that is, the reverse H\"{o}lder inequality.

A function $f$ is said to satisfy the reverse H\"{o}lder class if there exists a constant $C>0$ such that for any ball $B\subset\mathbb{R}^{n}$,
\begin{equation} \label{1.1}
\sup\limits_{B\ni x}|f(x)-f_{B}|\leq\frac{C}{|B|}\int_{B}|f(x)-f_{B}|dx.
\end{equation}
In fact, the reverse H\"{o}lder class contains many kinds functions, such as polynomial functions \cite{FC} and harmonic functions \cite{GT}. For more information about the reverse H\"{o}lder class, see also \cite{CN,HSV} for example.
\begin{theorem}[\cite{SL3}]
Let $1<p<\infty$, $-\frac{1}{p}\leq\lambda<0$, $-\frac{1}{p_{j}}\leq\lambda_{j}<0, j=1,2$, $\frac{1}{p}=\sum^{2}_{j=1}\frac{1}{p_{j}}$, $\lambda=\sum^{2}_{j=1}\lambda_{j}$ and let $b$ satisfy (\ref{1.1}). Then the following statements are equivalent:

(a)\quad $b\in\mathcal{C}^{p_{1},\lambda_{1}}(\mathbb{R}^{n})$;

(b)\quad Both $\mathcal{H}_{b}$ and $\mathcal{H}^{*}_{b}$ are bounded operators from $\mathcal{B}^{p_{2},\lambda_{2}}(\mathbb{R}^n)$ to $\mathcal{B}^{p,\lambda}(\mathbb{R}^n)$.
\end{theorem}
In addition, there are many rich and significant results about the characterizations of $\mathcal{C}^{p,\lambda}(\mathbb{R}^{n})$ via the boundedness of the Calder\'{o}n-Zygmund singular operator, fractional integral and Hardy operator, see e.g \cite{LS,SL1,SL2}.

In recent years, the topic of function spaces with the mixed norm has also increasing considerable attention and much developments, which stemmed from the fact that the more accurate structure of the mixed-norm function spaces than the corresponding classical function spaces such that the mixed-norm function spaces have more extensive applications in various areas such as potential analysis, harmonic analysis and PDEs, see \cite{AI,DKP,KC,KD,KN,LLW,LW,LTW1,LTW2} and the references therein. Next, we first recall the definitions of some function spaces with the mixed norm, and then state the motivation and main work of this paper.

In 1961, Benedek and Panzone \cite{BP} first introduced the following mixed Lebesgue spaces. For $0<\vec{p}=(p_{1},p_{2},\cdots,p_{n})<\infty$ (that is, $0<p_{i}<\infty$ for $i=1,2,\cdots,n$), the mixed Lebesgue norm $\|\cdot\|_{L^{\vec{p}}(\mathbb{R}^{n})}$ is defined by
$$\|f\|_{L^{\vec{p}}(\mathbb{R}^{n})}:=\Bigg(\int_{\mathbb{R}}\cdots\bigg(\int_{\mathbb{R}}\big(\int_{\mathbb{R}}|f(x_{1},x_{2},\cdots,x_{n})|^{p_{1}}
dx_{1}\big)^{\frac{p_{2}}{p_{1}}}dx_{2}\bigg)^{\frac{p_{3}}{p_{2}}}\cdots dx_{n}\Bigg)^{\frac{1}{p_{n}}},$$
where $f:\mathbb{R}^{n}\rightarrow\mathbb{C}$ is a measurable function. If $p_{i}=\infty$, for some $i=1,\cdots,n$, we only need to make some suitable modifications. The mixed Lebesgue space $L^{\vec{p}}(\mathbb{R}^{n})$ is defined as the set of all $f\in\mathscr{M}(\mathbb{R}^{n})$(the class of Lebesgue measurable functions) with $\|f\|_{L^{\vec{p}}(\mathbb{R}^{n})}<\infty$. Note that, if we take $p_{1}=p_{2}=\cdots=p_{n}=p$, then the mixed Lebesgue space $L^{\vec{p}}(\mathbb{R}^{n})$ reduces to the classical Lebesgue space
$L^{p}(\mathbb{R}^{n})$.

In 2019, Nagayama replaced the $L^{p}(\mathbb{R}^{n})$ norm by the mixed $L^{\vec{p}}(\mathbb{R}^{n})$ norm in the definition of classical Morrey space $M^{q}_{p}(\mathbb{R}^{n})$ in \cite{M1}, introduced the mixed Morrey space $M^{q}_{\vec{p}}(\mathbb{R}^{n})$ in \cite{NT1}. Recently, the authors \cite{LZ1} introduced mixed $\lambda$-central Morrey space $\mathcal{B}^{\vec{p},\lambda}(\mathbb{R}^n)$ and mixed $\lambda$-central bounded mean oscillation space $\mathrm{CMO}^{\vec{p},\lambda}(\mathbb{R}^{n})$ as follows.

\begin{definition}[\cite{LZ1}]
Let $1<\vec{p}<\infty$ and $\lambda\in\mathbb{R}$,
the mixed $\lambda$-central Morrey space $\mathcal{B}^{\vec{p},\lambda}(\mathbb{R}^n)$ is defined by
$$\mathcal{B}^{\vec{p},\lambda}(\mathbb{R}^n):=\{f\in \mathscr{M}(\mathbb{R}^{n}):\|f\|_{\mathcal{B}^{\vec{p},\lambda}(\mathbb{R}^n)}<\infty\},$$
where
$$\|f\|_{\mathcal{B}^{\vec{p},\lambda}(\mathbb{R}^n)}:=\sup\limits_{r>0}\frac{\|f\chi_{B(0,r)}\|_{L^{\vec{p}}(\mathbb{R}^{n})}}{|B(0,r)|^{\lambda}\|\chi_{B(0,r)}\|_{L^{\vec{p}}(\mathbb{R}^{n})}}.$$
\end{definition}

\begin{definition}[\cite{LZ1}]
Let $1<\vec{p}<\infty$ and $\lambda<\frac{1}{n}$,
the mixed $\lambda$-central bounded mean oscillation space $\mathrm{CMO}^{\vec{p},\lambda}(\mathbb{R}^{n})$ is defined by
$$\mathrm{CMO}^{\vec{p},\lambda}(\mathbb{R}^{n}):=\{f\in \mathscr{M}(\mathbb{R}^{n}):\|f\|_{\mathrm{CMO}^{\vec{p},\lambda}(\mathbb{R}^{n})}<\infty\},$$
where
$$\|f\|_{\mathrm{CMO}^{\vec{p},\lambda}(\mathbb{R}^{n})}:=
\sup\limits_{r>0}\frac{\|(f-f_{B(0,r)})\chi_{B(0,r)}\|_{L^{\vec{p}}(\mathbb{R}^{n})}}{|B(0,r)|^{\lambda}\|\chi_{B(0,r)}\|_{L^{\vec{p}}(\mathbb{R}^{n})}}.$$
\end{definition}

In \cite{LZ1}, the authors established the boundedness of some classical operators on the mixed $\lambda$-central Morrey space $\mathcal{B}^{\vec{p},\lambda}(\mathbb{R}^{n})$, such as  Calder\'{o}n-Zygmund singular integral operator, fractional integral operator and their commutators. Subsequently, the authors further extended  to the case of multi-sublinear operators and commutators in \cite{LZ2}. Furthermore, a characterization of mixed $\lambda$-central bounded mean oscillation space $\mathrm{CMO}^{\vec{p},\lambda}(\mathbb{R}^{n})$ was obtained in  \cite{LZ3} via the boundedness of the commutators $\mathcal{H}_{b}$ and $\mathcal{H}^{*}_{b}$ on mixed Lebesgue spaces. Meanwhile, the boundedness of commutators $\mathcal{H}_{b}$ and $\mathcal{H}^{*}_{b}$ on mixed $\lambda$-central Morrey space $\mathcal{B}^{\vec{p},\lambda}(\mathbb{R}^{n})$ was investigated.
\begin{theorem}[\cite{LZ3}]
Let $1<\vec{p},\vec{q}<\infty$, $0\leq\lambda<\frac{1}{n}$ satisfy condition $\lambda=\frac{1}{n}\sum^{n}_{i=1}\frac{1}{p_{i}}-\frac{1}{n}\sum^{n}_{i=1}\frac{1}{q_{i}}$. Then both $\mathcal{H}_{b}$ and $\mathcal{H}^{*}_{b}$ are bounded from $L^{\vec{p}}(\mathbb{R}^{n})$ to $L^{\vec{q}}(\mathbb{R}^{n})$ if and only if
$$b\in \mathrm{CMO}^{\vec{q},\lambda}(\mathbb{R}^{n})\cap\mathrm{CMO}^{\vec{p}',\lambda}(\mathbb{R}^{n}).$$
\end{theorem}

In view of this, the facts mentioned above pave the way for further research. A question that arises naturally is whether or not the results of Theorem 1.2 are also valid on the corresponding mixed $\lambda$-central Morrey spaces.

The main purpose in this paper is to solve the above problem. To achieve this goal, we first rewrite the central Campanato space as new function space with the mixed norm.

The mixed central Campanato space $\mathfrak{C}^{\vec{p},\lambda}(\mathbb{R}^{n})$ is defined as follows.

\begin{definition}
Let $1\leq \vec{p}<\infty$ and $-\frac{1}{n}\sum^{n}_{i=1}\frac{1}{p_{i}}<\lambda<\frac{1}{n}$,
we define the mixed central Campanato space $\mathfrak{C}^{\vec{p},\lambda}(\mathbb{R}^{n})$ as follows:
$$\mathfrak{C}^{\vec{p},\lambda}(\mathbb{R}^{n}):=\{f\in\mathscr{M}(\mathbb{R}^{n}):\|f\|_{\mathfrak{C}^{\vec{p},\lambda}(\mathbb{R}^{n})}<\infty\},$$
where
$$\|f\|_{\mathfrak{C}^{\vec{p},\lambda}(\mathbb{R}^{n})}:=
\sup\limits_{r>0}\frac{\|(f-f_{B(0,r)})\chi_{B(0,r)}\|_{L^{\vec{p}}(\mathbb{R}^{n})}}{|B(0,r)|^{\lambda}\|\chi_{B(0,r)}\|_{L^{\vec{p}}(\mathbb{R}^{n})}}.$$
\end{definition}

\hspace{-22pt} {\bf Remark 1.1} {(i) If $\vec{p}=p$, that is, $p_{i}=p, i=1,2,\cdots,n$, then
$$\mathfrak{C}^{\vec{p},\lambda}(\mathbb{R}^{n})=CMO^{p,\lambda}(\mathbb{R}^{n}).$$
(ii) For any $1\leq\vec{p}_1<\vec{p}_2<\infty$,  that is, $p_{1i}<p_{2i}, i=1,2,\cdots,n$, then
$$\mathfrak{C}^{\vec{p}_{2},\lambda}(\mathbb{R}^{n})\subset\mathfrak{C}^{\vec{p}_{1},\lambda}(\mathbb{R}^{n}).$$
(iii) For any $1\leq\vec{p}_1, \vec{p}_2<\infty$, if exist some $i_0$ such that $p_{1i_{0}}<p_{2i_{0}}$ and $p_{2i}<p_{1i}$ $i=1,\cdots,i_{0}-1,i_{0}+1,\cdots,n$, then the following relationship
$$\mathfrak{C}^{\vec{p}_{2},\lambda}(\mathbb{R}^{n})\subset\mathfrak{C}^{\vec{p}_{1},\lambda}(\mathbb{R}^{n})$$
does not hold.\\
(iv)  For any $1\leq\vec{p}<\infty$, if $-\frac{1}{n}\sum^{n}_{i=1}\frac{1}{p_{i}}<\lambda_{1}<\lambda_{2}<\frac{1}{n}$, then
$$\mathfrak{C}^{\vec{p},\lambda_{1}}(\mathbb{R}^{n})\subset\mathfrak{C}^{\vec{p},\lambda_{2}}(\mathbb{R}^{n}).$$
(v) For any $1\leq\vec{p}<\infty$, if $\lambda=0$, then the mixed central Campanato space
$\mathfrak{C}^{\vec{p},0}(\mathbb{R}^{n})$ reduces to the mixed central bounded mean oscillation space $CMO^{\vec{p}}(\mathbb{R}^{n})$, which was defined by Wei in \cite{WM2}. Also, similar $BMO(\mathbb{R}^{n})$, the following equivalent norm of $CMO^{\vec{p}}(\mathbb{R}^{n})$
$$\|f\|_{\mathrm{CMO}^{\vec{p}}(\mathbb{R}^{n})}:=\sup\limits_{r>0}
\inf\limits_{c\in\mathbb{C}}\frac{\|(f-c)\chi_{B(0,r)}\|_{L^{\vec{p}}(\mathbb{R}^{n})}}{\|\chi_{B(0,r)}\|_{L^{\vec{p}}(\mathbb{R}^{n})}}$$
is valid.
}

In order to establish some characterizations of the mixed central Campanato space $\mathfrak{C}^{\vec{p},\lambda}(\mathbb{R}^{n})$ via the boundedness of the commutators of Hardy type on the mixed $\lambda$-central Morrey space $\mathcal{B}^{\vec{p},\lambda}(\mathbb{R}^n)$, we first need to modify the reverse H\"{o}lder class as follows: Consider functions $f$ such that there exists a constant $C>0$ such that for any $B\subset\mathbb{R}^{n}$,
\begin{equation} \label{1.2}
\sup\limits_{B\ni x}|f(x)-f_{B}|\leq C\frac{\|(f-f_{B})\chi_{B}\|_{L^{1}(\mathbb{R}^{n})}}{\|\chi_{B}\|_{L^{1}(\mathbb{R}^{n})}},
\end{equation}
here $\|\cdot\|_{L^{1}(\mathbb{R}^{n})}$ denotes the mixed $L^{1}$ norm. A light abuse of terminology, we can also call it the reverse H\"{o}lder class.

The main results in this paper can be formulated as follows.
\begin{theorem}
Let $1<\vec{p}, \vec{p_{2}}<\infty$, $-\frac{1}{n}\sum^{n}_{i=1}\frac{1}{p_{i}}<\lambda<0$, $-\frac{1}{n}\sum^{n}_{i=1}\frac{1}{p_{ji}}<\lambda_{j}<0$, $j=1,2$, $\frac{1}{\vec{p}}=\frac{1}{\vec{p_{1}}}+\frac{1}{\vec{p_{2}}}$, $\lambda=\lambda_{1}+\lambda_{2}$ and suppose $b$ satisfies (\ref{1.2}). Then the following statements are equivalent:

(a)\quad $b\in\mathfrak{C}^{\vec{p_{1}},\lambda_{1}}(\mathbb{R}^{n})$;

(b)\quad Both $\mathcal{H}_{b}$ and $\mathcal{H}^{*}_{b}$ are bounded operators from $\mathcal{B}^{\vec{p_{2}},\lambda_{2}}(\mathbb{R}^n)$ to $\mathcal{B}^{\vec{p},\lambda}(\mathbb{R}^n)$.
\end{theorem}

Under some stronger conditions on $\lambda$ and $\vec{p}$, the following result can be obtained if remove the assumption that $b$ satisfies the condition (\ref{1.2}).

\begin{theorem}
Let $2<\vec{p}<\infty$, $-\frac{1}{n}\sum^{n}_{i=1}\frac{1}{2p_{i}}<\lambda<0$. Then the following statements are equivalent:

(a)\quad $b\in\mathfrak{C}^{\vec{p},\lambda}(\mathbb{R}^{n})$;

(b)\quad Both $\mathcal{H}_{b}$ and $\mathcal{H}^{*}_{b}$ are bounded operators from $\mathcal{B}^{\vec{p},\lambda}(\mathbb{R}^n)$ to $\mathcal{B}^{\vec{p},2\lambda}(\mathbb{R}^n)$.
\end{theorem}

In 2007, Fu et al. \cite{FLW} gave definitions of the $n$-dimensional fractional Hardy operator $\mathcal{H}_{\alpha}$, its dual operator $\mathcal{H}^{*}_{\alpha}$ and their commutators. Let $f$ be a locally integrable function in $\mathbb{R}^{n}$ and $0<\alpha<n$. The $n$-dimensional fractional Hardy operator can be defined by
$$\mathcal{H}_{\alpha}f(x)=\frac{1}{|x|^{n-\alpha}}\int_{|t|<|x|}f(t)dt, \quad x\in\mathbb{R}^{n}\backslash\{0\}.$$
The dual operator of $\mathcal{H}_{\alpha}$ is $\mathcal{H}^{*}_{\alpha}$,
$$\mathcal{H}^{*}_{\alpha}f(x)=\int_{|t|\geq|x|}\frac{f(t)}{|t|^{n-\alpha}}dt, \quad x\in\mathbb{R}^{n}\backslash\{0\}.$$
The commutators of $\mathcal{H}_{\alpha}$ and $\mathcal{H}^{*}_{\alpha}$ are defined by
$$\mathcal{H}_{\alpha,b}f=b\mathcal{H}_{\alpha}f-\mathcal{H}_{\alpha}(fb),\quad
\mathcal{H}^{*}_{\alpha,b}f=b\mathcal{H}^{*}_{\alpha}f-\mathcal{H}^{*}_{\alpha}(fb).$$
Furthermore, they established some characterizations of $\mathrm{CMO}^{p}(\mathbb{R}^{n})(1<p<\infty)$ via the boundedness of $\mathcal{H}_{\alpha,b}$ on both Lebesgue spaces and Herz spaces by the dual technology.

Next, we establish some characterizations of $\mathfrak{C}^{\vec{p},\lambda}(\mathbb{R}^{n})$ with $\lambda<0$ via the boundedness of $\mathcal{H}_{\alpha,b}$ and $\mathcal{H}^{*}_{\alpha,b}$.

\begin{theorem}
Let $\vec{p}, \lambda, \vec{p_{j}}, \lambda_{j}, j=1,2$, $b$ be as in Theorem 1.4, $0<\alpha<\min\{n(1-\frac{1}{n}\sum^{n}_{i=1}\frac{1}{p_{i}}), n(\lambda_{2}+\frac{1}{n}\sum^{n}_{i=1}\frac{1}{p_{2i}})\}$, and $\beta=\lambda_{2}-\frac{\alpha}{n}$. Then the following statements are equivalent:

(a)\quad $b\in\mathfrak{C}^{\vec{p_{1}},\lambda_{1}}(\mathbb{R}^{n})$;

(b)\quad Both $\mathcal{H}_{\alpha,b}$ and $\mathcal{H}^{*}_{\alpha,b}$ are bounded operators from $\mathcal{B}^{\vec{p_{2}},\beta}(\mathbb{R}^n)$ to $\mathcal{B}^{\vec{p},\lambda}(\mathbb{R}^n)$.
\end{theorem}

\begin{theorem}
Let $2<\vec{p}<\infty$, $-\frac{1}{n}\sum^{n}_{i=1}\frac{1}{2p_{i}}<\lambda<0$, $0<\alpha<\min\{n(1-\frac{1}{n}\sum^{n}_{i=1}\frac{1}{p_{i}}), n(\lambda+\frac{1}{n}\sum^{n}_{i=1}\frac{1}{p_{i}})\}$, and $\beta=\lambda-\frac{\alpha}{n}$. Then the following statements are equivalent:

(a)\quad $b\in\mathfrak{C}^{\vec{p},\lambda}(\mathbb{R}^{n})$;

(b)\quad Both $\mathcal{H}_{\alpha,b}$ and $\mathcal{H}^{*}_{\alpha,b}$ are bounded operators from $\mathcal{B}^{\vec{p},\beta}(\mathbb{R}^n)$ to $\mathcal{B}^{\vec{p},2\lambda}(\mathbb{R}^n)$.
\end{theorem}

The rest of this paper is organized as follows. In Section 2, we will prove Theorems 1.4-1.5. Some characterizations of the mixed central bounded mean oscillation space $CMO^{\vec{p}}(\mathbb{R}^n)$ and the boundedness of Hardy type operators on mixed $\lambda$-central Morrey space $\mathcal{B}^{\vec{p},\lambda}(\mathbb{R}^n)$ are obtained in Section 3.

As a rule, for $1\leq\vec{p}\leq\infty$, let $$\frac{1}{\vec{p}}=(\frac{1}{p_{1}},\frac{1}{p_{2}},\cdots,\frac{1}{p_{n}}), \quad\vec{p}'=(p_{1}',p_{2}',\cdots,p_{n}'),$$
where $p_{i}'=\frac{p_{i}}{p_{i}-1} (i=1,2,\cdots,n)$ is the conjugate exponent of $p_i$. For any set $E\subset\mathbb{R}^n$, ${\chi}_{E}$ denotes its characteristic function and $E^{c}$ denotes its complementary set, we also denote the Lebesgue measure by $|E|$. Throughout this paper, We use the symbol $f\lesssim g$ to denote there exists a positive constant $C$ such that $f\leq Cg $, and the notation $f\thickapprox g$ means that there exist positive constants $C_1, C_2$ such that $C_1 g\leq f\leq C_2g$.

\section{Proofs of the main results}

In this section, we provide the proofs of  main results. In fact, for Theorem 1.4, it is enough to grasp the core techniques adopted in the proof, and Theorem 1.5 is an expected conclusion without the assumption of the reverse H\"{o}lder class. Therefore, we only need to present the proofs of Theorems 1.4-1.5, the other proofs require only a simple modification. In what follows, for simplicity, we write
$$B_{k}=\{x\in\mathbb{R}^{n}: |x|\leq 2^{k}\}, \quad C_{k}=B_{k}\backslash B_{k-1}, \quad \chi_{k}=\chi_{C_{k}},\quad k\in\mathbb{Z},$$
where $\chi_{C_{k}}$ is the characteristic function of a set $C_{k}$.

\begin{proof}[\bf Proof of Theorem 1.4]

The proof of this theorem can be divided into two steps.

$(a)\Rightarrow(b)$ Let $f$ be a function in $\mathcal{B}^{\vec{p_{2}},\lambda_{2}}(\mathbb{R}^n)$, $b\in\mathfrak{C}^{\vec{p_{1}},\lambda_{1}}(\mathbb{R}^{n})$. For a fixed ball $B=B(0,r)\subset\mathbb{R}^{n}$, there is no loss of generality in assuming $B(0,r)=B_{k_{0}}$ with $k_{0}\in\mathbb{Z}$, We just need to show that the facts
\begin{equation} \label{2.1}
\|(\mathcal{H}_{b}f)\chi_{B_{k_{0}}}\|_{\vec{p}}\lesssim|B_{k_{0}}|^{\lambda}\|\chi_{B_{k_{0}}}\|_{\vec{p}}\|b\|_{\mathfrak{C}^{\vec{p_{1}},\lambda_{1}}(\mathbb{R}^{n})}
\|f\|_{\mathcal{B}^{\vec{p_{2}},\lambda_{2}}(\mathbb{R}^{n})}
\end{equation}
and
\begin{equation} \label{2.2}
\|(\mathcal{H}^{*}_{b}f)\chi_{B_{k_{0}}}\|_{\vec{p}}\lesssim|B_{k_{0}}|^{\lambda}\|\chi_{B_{k_{0}}}\|_{\vec{p}}\|b\|_{\mathfrak{C}^{\vec{p_{1}},\lambda_{1}}(\mathbb{R}^{n})}
\|f\|_{\mathcal{B}^{\vec{p_{2}},\lambda_{2}}(\mathbb{R}^{n})}
\end{equation}
hold.

The definition of $\mathcal{H}_{b}$ and the Minkowski inequality on mixed Lebesgue space (see \cite{BP}) give that
\begin{align*}
\|\mathcal{H}_{b}f(\cdot)\chi_{B_{k_{0}}}(\cdot)\|_{\vec{p}}&=
\Big\|\chi_{B_{k_{0}}}(\cdot)\frac{1}{|\cdot|^{n}}\int_{|y|<|\cdot|}\big(b(\cdot)-b(y)\big)f(y)dy\Big\|_{\vec{p}}\\
&\leq\sum^{k_{0}}_{k=-\infty}\Big\|\chi_{C_{k}}(\cdot)\frac{1}{|\cdot|^{n}}\int_{B_{k}}\big(b(\cdot)-b(y)\big)f(y)dy\Big\|_{\vec{p}}
\end{align*}
\begin{align*}
&\lesssim\sum^{k_{0}}_{k=-\infty}\Big\|\chi_{C_{k}}(\cdot)\frac{1}{|\cdot|^{n}}\sum^{k}_{i=-\infty}
\int_{C_{i}}\big(b(\cdot)-b_{B_{k}}\big)f(y)dy\Big\|_{\vec{p}}\\
&\quad+\sum^{k_{0}}_{k=-\infty}\Big\|\chi_{C_{k}}(\cdot)\frac{1}{|\cdot|^{n}}\sum^{k}_{i=-\infty}
\int_{C_{i}}\big(b(y)-b_{B_{k}}\big)f(y)dy\Big\|_{\vec{p}}\\
&=:I+II.
\end{align*}

First, we estimate $I$. For $\frac{1}{\vec{p}}=\frac{1}{\vec{p_{1}}}+\frac{1}{\vec{p_{2}}}$ and $1=\frac{1}{\vec{p_{2}}}+\frac{1}{\vec{p_{2}}'}$, by the H\"{o}lder inequality on mixed Lebesgue space (see \cite{BP}), we get
\begin{align*}
I&\lesssim\sum^{k_{0}}_{k=-\infty}2^{-kn}\Big\|\chi_{C_{k}}\big(b-b_{B_{k}}\big)\sum^{k}_{i=-\infty}\int_{C_{i}}f(y)dy\Big\|_{\vec{p}}\\
&\leq\sum^{k_{0}}_{k=-\infty}2^{-kn}\Big\|\chi_{C_{k}}\big(b-b_{B_{k}}\big)\Big\|_{\vec{p_{1}}}
\sum^{k}_{i=-\infty}\Big\|\chi_{C_{k}}\int_{C_{i}}f(y)dy\Big\|_{\vec{p_{2}}}\\
&\leq\sum^{k_{0}}_{k=-\infty}2^{-kn}\big\|\chi_{B_{k}}\big(b-b_{B_{k}}\big)\big\|_{\vec{p_{1}}}\sum^{k}_{i=-\infty}
\big\|\chi_{B_{k}}\big\|_{\vec{p_{2}}}\big\|f\chi_{B_{i}}\big\|_{\vec{p_{2}}}\big\|\chi_{B_{i}}\big\|_{\vec{p_{2}}'}\\
&\leq\|b\|_{\mathfrak{C}^{\vec{p_{1}},\lambda_{1}}(\mathbb{R}^{n})}\|f\|_{\mathcal{B}^{\vec{p_{2}},\lambda_{2}}(\mathbb{R}^{n})}
\sum^{k_{0}}_{k=-\infty}2^{-kn}|B_{k}|^{\lambda_{1}}\big\|\chi_{B_{k}}\big\|_{\vec{p_{1}}}\\
&\quad\times\sum^{k}_{i=-\infty}|B_{i}|^{\lambda_{2}}
\big\|\chi_{B_{i}}\big\|_{\vec{p_{2}}}\big\|\chi_{B_{i}}\big\|_{\vec{p_{2}}'}\big\|\chi_{B_{k}}\big\|_{\vec{p_{2}}}\\
&\lesssim\|b\|_{\mathfrak{C}^{\vec{p_{1}},\lambda_{1}}(\mathbb{R}^{n})}\|f\|_{\mathcal{B}^{\vec{p_{2}},\lambda_{2}}(\mathbb{R}^{n})}
\sum^{k_{0}}_{k=-\infty}2^{-kn}|B_{k}|^{\lambda_{1}}\|\chi_{B_{k}}\|_{\vec{p}}\sum^{k}_{i=-\infty}|B_{i}|^{\lambda_{2}+1}\\
&\lesssim|B_{k_{0}}|^{\lambda}\|\chi_{B_{k_{0}}}\|_{\vec{p}}
\|b\|_{\mathfrak{C}^{\vec{p_{1}},\lambda_{1}}(\mathbb{R}^{n})}\|f\|_{\mathcal{B}^{\vec{p_{2}},\lambda_{2}}(\mathbb{R}^{n})},
\end{align*}
where
$$\|\chi_{B_{k}}\|_{\vec{p_{1}}}\|\chi_{B_{k}}\|_{\vec{p_{2}}}\approx|B_{k}|^{\frac{1}{n}\sum^{n}_{i=1}\frac{1}{p_{1i}}+
\frac{1}{n}\sum^{n}_{i=1}\frac{1}{p_{2i}}}\approx\|\chi_{B_{k}}\|_{\vec{p}},$$
and we have used the conditions $-\frac{1}{n}\sum^{n}_{i=1}\frac{1}{p_{i}}<\lambda<0$, $-\frac{1}{n}\sum^{n}_{i=1}\frac{1}{p_{2i}}<\lambda_{2}<0$ and $\lambda=\lambda_{1}+\lambda_{2}$.

Next, by H\"{o}lder's inequality on mixed Lebesgue space, the fact $1=\frac{1}{\vec{p}}+\frac{1}{\vec{p}'}$ and $\frac{1}{\vec{p}}=\frac{1}{\vec{p_{1}}}+\frac{1}{\vec{p_{2}}}$ allow us to estimate the term $II$ as
\begin{align*}
II&\lesssim\sum^{k_{0}}_{k=-\infty}2^{-kn}\Big\|\chi_{C_{k}}\sum^{k}_{i=-\infty}\int_{B_{i}}\big(b(y)-b_{B_{k}}\big)f(y)dy\Big\|_{\vec{p}}\\
&\leq\sum^{k_{0}}_{k=-\infty}2^{-kn}\Big\|\chi_{C_{k}}
\sum^{k}_{i=-\infty}\big\|(b-b_{B_{k}})f\chi_{B_{i}}\big\|_{\vec{p}}\big\|\chi_{B_{i}}\big\|_{\vec{p}'}\Big\|_{\vec{p}}\\
&\leq\sum^{k_{0}}_{k=-\infty}2^{-kn}\big\|\chi_{B_{k}}\big\|_{\vec{p}}
\sum^{k}_{i=-\infty}\big\|(b-b_{B_{k}})\chi_{B_{i}}\big\|_{\vec{p_{1}}}\big\|f\chi_{B_{i}}\big\|_{\vec{p_{2}}}\big\|\chi_{B_{i}}\big\|_{\vec{p}'}\\
&\leq\|b\|_{\mathfrak{C}^{\vec{p_{1}},\lambda_{1}}(\mathbb{R}^{n})}\|f\|_{\mathcal{B}^{\vec{p_{2}},\lambda_{2}}(\mathbb{R}^{n})}
\sum^{k_{0}}_{k=-\infty}2^{-kn}|B_{k}|^{\lambda_{1}}\big\|\chi_{B_{k}}\big\|_{\vec{p_{1}}}\big\|\chi_{B_{k}}\big\|_{\vec{p}}
\end{align*}
\begin{align*}
&\quad\times\sum^{k}_{i=-\infty}|B_{i}|^{\lambda_{2}}\big\|\chi_{B_{i}}\big\|_{\vec{p_{2}}}\big\|\chi_{B_{i}}\big\|_{\vec{p}'}\\
&\lesssim\|b\|_{\mathfrak{C}^{\vec{p_{1}},\lambda_{1}}(\mathbb{R}^{n})}\|f\|_{\mathcal{B}^{\vec{p_{2}},\lambda_{2}}(\mathbb{R}^{n})}
\sum^{k_{0}}_{k=-\infty}2^{kn(\lambda_{1}+\frac{1}{n}\sum^{n}_{i=1}\frac{1}{p_{i}}+\frac{1}{n}\sum^{n}_{i=1}\frac{1}{p_{1i}}-1)}\\
&\quad\times\sum^{k}_{i=-\infty}2^{in(\lambda_{2}+\frac{1}{n}\sum^{n}_{i=1}\frac{1}{p_{i}'}+\frac{1}{n}\sum^{n}_{i=1}\frac{1}{p_{2i}})}\\
&\leq|B_{k_{0}}|^{\lambda}\|\chi_{B_{k_{0}}}\|_{\vec{p}}
\|b\|_{\mathfrak{C}^{\vec{p_{1}},\lambda_{1}}(\mathbb{R}^{n})}\|f\|_{\mathcal{B}^{\vec{p_{2}},\lambda_{2}}(\mathbb{R}^{n})},
\end{align*}
and we have used the fact that $-\frac{1}{n}\sum^{n}_{i=1}\frac{1}{p_{2i}}<\lambda_{2}<0$, $\lambda=\lambda_{1}+\lambda_{2}$, and $-\frac{1}{n}\sum^{n}_{i=1}\frac{1}{p_{i}}<\lambda<0$.

Combining the estimates of $I$ and $II$, (\ref{2.1}) is obtained.

Now, let's prove (\ref{2.2}) in this position. Note that
\begin{align*}
\|\mathcal{H}^{*}_{b}f(\cdot)\chi_{B_{k_{0}}}(\cdot)\|_{\vec{p}}&=
\Big\|\chi_{B_{k_{0}}}(\cdot)\int_{|y|\geq|\cdot|}\frac{\big(b(\cdot)-b(y)\big)f(y)}{|y|^{n}}dy\Big\|_{\vec{p}}\\
&\leq\Big\|\chi_{B_{k_{0}}}(\cdot)\int_{2^{k_{0}n}\geq|y|\geq|\cdot|}\frac{|b(\cdot)-b(y)||f(y)|}{|y|^{n}}dy\Big\|_{\vec{p}}\\
&\quad+\Big\|\chi_{B_{k_{0}}}(\cdot)\int_{|y|>2^{k_{0}n}}\frac{|b(\cdot)-b(y)||f(y)|}{|y|^{n}}dy\Big\|_{\vec{p}}\\
&=:M+N.
\end{align*}
The term $M$ can be handled in much the same way as that of (\ref{2.1}), the only difference being in the analysis of the term $N$. Analysis similar to that of $\mathcal{H}_{b}$ shows
\begin{align*}
M&\leq\Big\|\chi_{B_{k_{0}}}(\cdot)\frac{1}{|\cdot|^{n}}\int_{|y|\leq2^{k_{0}n}}|b(\cdot)-b(y)||f(y)|dy\Big\|_{\vec{p}}\\
&\leq\sum^{k_{0}}_{k=-\infty}\Big\|\chi_{C_{k}}(\cdot)\frac{1}{|\cdot|^{n}}\sum^{k}_{i=-\infty}
\int_{C_{i}}|\big(b(\cdot)-b(y)\big)f(y)|dy\Big\|_{\vec{p}}\\
&\lesssim\sum^{k_{0}}_{k=-\infty}2^{-kn}\Big\|\chi_{B_{k}}(\cdot)\sum^{k}_{i=-\infty}\int_{B_{i}}|\big(b(\cdot)-b(y)\big)f(y)|dy\Big\|_{\vec{p}}\\
&\lesssim|B_{k_{0}}|^{\lambda}\|\chi_{B_{k_{0}}}\|_{\vec{p}}
\|b\|_{\mathfrak{C}^{\vec{p_{1}},\lambda_{1}}(\mathbb{R}^{n})}\|f\|_{\mathcal{B}^{\vec{p_{2}},\lambda_{2}}(\mathbb{R}^{n})}.
\end{align*}
For the term $N$, using the Minkowski inequality on mixed Lebesgue space, we write
\begin{align*}
N&\leq\Big\|\chi_{B_{k_{0}}}(\cdot)\sum^{\infty}_{k=k_{0}}\int_{C_{k}}\frac{|b(\cdot)-b_{B_{k_{0}}}||f(y)|}{|y|^{n}}dy\Big\|_{\vec{p}}\\
&\quad+\Big\|\chi_{B_{k_{0}}}(\cdot)\sum^{\infty}_{k=k_{0}}\int_{C_{k}}\frac{|b(y)-b_{B_{k_{0}}}||f(y)|}{|y|^{n}}dy\Big\|_{\vec{p}}\\
&=:N_{1}+N_{2}.
\end{align*}
For $\frac{1}{\vec{p}}=\frac{1}{\vec{p_{1}}}+\frac{1}{\vec{p_{2}}}$ and $1=\frac{1}{\vec{p_{2}}}+\frac{1}{\vec{p_{2}}'}$, by the H\"{o}lder inequality on mixed Lebesgue space, the fact $-\frac{1}{n}\sum^{n}_{i=1}\frac{1}{p_{2i}}<\lambda_{2}<0$ and $\lambda=\lambda_{1}+\lambda_{2}$ deduce that

\begin{align*}
N_{1}&=\Big\|\chi_{B_{k_{0}}}(\cdot)\big(b(\cdot)-b_{B_{k_{0}}}\big)\sum^{\infty}_{k=k_{0}}\int_{C_{k}}\frac{|f(y)|}{|y|^{n}}dy\Big\|_{\vec{p}}\\
&\leq\Big\|\chi_{B_{k_{0}}}\big(b-b_{B_{k_{0}}}\big)\Big\|_{\vec{p_{1}}}\Big\|\chi_{B_{k_{0}}}\Big\|_{\vec{p_{2}}}
\sum^{\infty}_{k=k_{0}}\Big\|\chi_{C_{k}}(\cdot)\frac{|f(\cdot)|}{|\cdot|^{n}}\Big\|_{\vec{p_{2}}}\Big\|\chi_{C_{k}}\Big\|_{\vec{p_{2}}'}\\
&\lesssim|B_{k_{0}}|^{\lambda_{1}}\|\chi_{B_{k_{0}}}\|_{\vec{p}}\|b\|_{\mathfrak{C}^{\vec{p_{1}},\lambda_{1}}(\mathbb{R}^{n})}
\sum^{\infty}_{k=k_{0}}2^{-kn}\Big\|f\chi_{B_{k}}\Big\|_{\vec{p_{2}}}\Big\|\chi_{B_{k}}\Big\|_{\vec{p_{2}}'}\\
&\lesssim|B_{k_{0}}|^{\lambda_{1}}\|\chi_{B_{k_{0}}}\|_{\vec{p}}\|b\|_{\mathfrak{C}^{\vec{p_{1}},\lambda_{1}}(\mathbb{R}^{n})}
\|f\|_{\mathcal{B}^{\vec{p_{2}},\lambda_{2}}(\mathbb{R}^{n})}\sum^{\infty}_{k=k_{0}}2^{-kn}|B_{k}|^{\lambda_{2}+1}\\
&\lesssim|B_{k_{0}}|^{\lambda}\|\chi_{B_{k_{0}}}\|_{\vec{p}}\|b\|_{\mathfrak{C}^{\vec{p_{1}},\lambda_{1}}(\mathbb{R}^{n})}
\|f\|_{\mathcal{B}^{\vec{p_{2}},\lambda_{2}}(\mathbb{R}^{n})}.
\end{align*}
To estimate $N_{2}$, we need the following decomposition
\begin{align*}
N_{2}&\leq\Big\|\chi_{B_{k_{0}}}(\cdot)\sum^{\infty}_{k=k_{0}}\int_{C_{k}}\frac{|b(y)-b_{B_{k}}||f(y)|}{|y|^{n}}dy\Big\|_{\vec{p}}
+\Big\|\chi_{B_{k_{0}}}(\cdot)\sum^{\infty}_{k=k_{0}}\int_{C_{k}}\frac{|b_{B_{k_{0}}}-b_{B_{k}}||f(y)|}{|y|^{n}}dy\Big\|_{\vec{p}}\\
&=:N_{21}+N_{22}.
\end{align*}
We first compute $N_{21}$. For $1=\frac{1}{\vec{p}}+\frac{1}{\vec{p}'}$ and $\frac{1}{\vec{p}}=\frac{1}{\vec{p_{1}}}+\frac{1}{\vec{p_{2}}}$, using the H\"{o}lder inequality on mixed Lebesgue space, we obtain
\begin{align*}
N_{21}&\lesssim\Big\|\chi_{B_{k_{0}}}(\cdot)\sum^{\infty}_{k=k_{0}}2^{-kn}\int_{B_{k}}|b(y)-b_{B_{k}}||f(y)|dy\Big\|_{\vec{p}}\\
&\leq\big\|\chi_{B_{k_{0}}}\big\|_{\vec{p}}\sum^{\infty}_{k=k_{0}}2^{-kn}
\big\|(b-b_{B_{k}})f\chi_{B_{k}}\big\|_{\vec{p}}\big\|\chi_{B_{k}}\big\|_{\vec{p}'}\\
&\leq\big\|\chi_{B_{k_{0}}}\big\|_{\vec{p}}\sum^{\infty}_{k=k_{0}}2^{-kn}
\big\|(b-b_{B_{k}})\chi_{B_{k}}\big\|_{\vec{p_{1}}}\big\|f\chi_{B_{k}}\big\|_{\vec{p_{2}}}\big\|\chi_{B_{k}}\big\|_{\vec{p}'}\\
&\lesssim\|\chi_{B_{k_{0}}}\|_{\vec{p}}\|b\|_{\mathfrak{C}^{\vec{p_{1}},\lambda_{1}}(\mathbb{R}^{n})}
\|f\|_{\mathcal{B}^{\vec{p_{2}},\lambda_{2}}(\mathbb{R}^{n})}\sum^{\infty}_{k=k_{0}}2^{-kn}|B_{k}|^{\lambda+1}\\
&\lesssim|B_{k_{0}}|^{\lambda}\|\chi_{B_{k_{0}}}\|_{\vec{p}}\|b\|_{\mathfrak{C}^{\vec{p_{1}},\lambda_{1}}(\mathbb{R}^{n})}
\|f\|_{\mathcal{B}^{\vec{p_{2}},\lambda_{2}}(\mathbb{R}^{n})},
\end{align*}
where we have used the conditions $\lambda=\lambda_{1}+\lambda_{2}$ and $-\frac{1}{n}\sum^{n}_{i=1}\frac{1}{p_{i}}<\lambda<0$.

For the term $N_{22}$, we claim first that for $k>k_{0}$,
$$|b_{B_{k_{0}}}-b_{B_{k}}|\lesssim2^{(k_{0}+1)n\lambda_{1}}\|b\|_{\mathfrak{C}^{\vec{p_{1}},\lambda_{1}}(\mathbb{R}^{n})}.$$
In fact, to do this, for $-\frac{1}{n}\sum^{n}_{i=1}\frac{1}{p_{1i}}<\lambda_{1}<0$, the H\"{o}lder inequality on mixed Lebesgue space and $1=\frac{1}{\vec{p_{1}}}+\frac{1}{\vec{p_{1}}'}$ show that
\begin{align*}
|b_{B_{k_{0}}}-b_{B_{k}}|&\leq\sum^{k-1}_{j=k_{0}}|b_{B_{j}}-b_{B_{j+1}}|
\end{align*}
\begin{align*}
&\leq\sum^{k-1}_{j=k_{0}}\frac{1}{|B_{j}|}\int_{B_{j+1}}|b(y)-b_{B_{j+1}}|dy\\
&\leq\sum^{k-1}_{j=k_{0}}\frac{1}{|B_{j}|}\big\|(b-b_{B_{j+1}})\chi_{B_{j+1}}\big\|_{\vec{p_{1}}}\big\|\chi_{B_{j+1}}\big\|_{\vec{p_{1}}'}\\
&\lesssim\|b\|_{\mathfrak{C}^{\vec{p_{1}},\lambda_{1}}(\mathbb{R}^{n})}\sum^{k-1}_{j=k_{0}}\frac{|B_{j+1}|^{\lambda_{1}+1}}{|B_{j}|}\\
&\lesssim\|b\|_{\mathfrak{C}^{\vec{p_{1}},\lambda_{1}}(\mathbb{R}^{n})}\sum^{k-1}_{j=k_{0}}2^{(j+1)n\lambda_{1}}\\
&\lesssim2^{(k_{0}+1)n\lambda_{1}}\|b\|_{\mathfrak{C}^{\vec{p_{1}},\lambda_{1}}(\mathbb{R}^{n})}.
\end{align*}
Therefore, for $-\frac{1}{n}\sum^{n}_{i=1}\frac{1}{p_{2i}}<\lambda_{2}<0$ and $\lambda=\lambda_{1}+\lambda_{2}$, we have
\begin{align*}
N_{22}&\lesssim\|b\|_{\mathfrak{C}^{\vec{p_{1}},\lambda_{1}}(\mathbb{R}^{n})}
\Big\|\chi_{B_{k_{0}}}(\cdot)\sum^{\infty}_{k=k_{0}}\int_{C_{k}}\frac{2^{k_{0}n\lambda_{1}}|f(y)|}{|y|^{n}}dy\Big\|_{\vec{p}}\\
&\lesssim\|b\|_{\mathfrak{C}^{\vec{p_{1}},\lambda_{1}}(\mathbb{R}^{n})}
\Big\|\chi_{B_{k_{0}}}(\cdot)2^{k_{0}n\lambda_{1}}\sum^{\infty}_{k=k_{0}}2^{-kn}\int_{B_{k}}|f(y)|dy\Big\|_{\vec{p}}\\
&\leq\|b\|_{\mathfrak{C}^{\vec{p_{1}},\lambda_{1}}(\mathbb{R}^{n})}|B_{k_{0}}|^{\lambda_{1}}\|\chi_{B_{k_{0}}}\|_{\vec{p}}
\sum^{\infty}_{k=k_{0}}2^{-kn}\|f\chi_{B_{k}}\|_{\vec{p_{2}}}\|\chi_{B_{k}}\|_{\vec{p_{2}'}}\\
&\lesssim|B_{k_{0}}|^{\lambda_{1}}\|\chi_{B_{k_{0}}}\|_{\vec{p}}\|b\|_{\mathfrak{C}^{\vec{p_{1}},\lambda_{1}}(\mathbb{R}^{n})}
\|f\|_{\mathcal{B}^{\vec{p_{2}},\lambda_{2}}(\mathbb{R}^{n})}\sum^{\infty}_{k=k_{0}}2^{kn\lambda_{2}}\\
&\leq|B_{k_{0}}|^{\lambda}\|\chi_{B_{k_{0}}}\|_{\vec{p}}\|b\|_{\mathfrak{C}^{\vec{p_{1}},\lambda_{1}}(\mathbb{R}^{n})}
\|f\|_{\mathcal{B}^{\vec{p_{2}},\lambda_{2}}(\mathbb{R}^{n})}.
\end{align*}
Summarizing, we get
$$N\lesssim|B_{k_{0}}|^{\lambda}\|\chi_{B_{k_{0}}}\|_{\vec{p}}\|b\|_{\mathfrak{C}^{\vec{p_{1}},\lambda_{1}}(\mathbb{R}^{n})}
\|f\|_{\mathcal{B}^{\vec{p_{2}},\lambda_{2}}(\mathbb{R}^{n})}.$$
Thus, $M$ together with $N$ implies that
$$\|(\mathcal{H}^{*}_{b}f)\chi_{B_{k_{0}}}\|_{\vec{p}}\lesssim|B_{k_{0}}|^{\lambda}\|\chi_{B_{k_{0}}}\|_{\vec{p}}\|b\|_{\mathfrak{C}^{\vec{p_{1}},\lambda_{1}}(\mathbb{R}^{n})}
\|f\|_{\mathcal{B}^{\vec{p_{2}},\lambda_{2}}(\mathbb{R}^{n})}.$$
This completes  (\ref{2.2}).

$(b)\Rightarrow(a)$ In this case, suppose that $\mathcal{H}_{b}$ and $\mathcal{H}^{*}_{b}$ are bounded operators from $\mathcal{B}^{\vec{p_{2}},\lambda_{2}}(\mathbb{R}^n)$ to $\mathcal{B}^{\vec{p},\lambda}(\mathbb{R}^n)$, we will show that
$$b\in\mathfrak{C}^{\vec{p_{1}},\lambda_{1}}(\mathbb{R}^{n}).$$
the key point of the proof is to construct a proper commutator. We are reduced to proving that for a fixed ball $B=B(0,r)$,
$$\frac{\|(b-b_{B})\chi_{B}\|_{\vec{p_{1}}}}{|B|^{\lambda_{1}}\|\chi_{B}\|_{\vec{p_{1}}}}\leq C.$$
We conclude from (\ref{1.2}) and the H\"{o}lder inequality on mixed Lebesgue space that
\begin{align*}
\frac{\|(b-b_{B})\chi_{B}\|_{\vec{p_{1}}}}{|B|^{\lambda_{1}}\|\chi_{B}\|_{\vec{p_{1}}}}&\leq\frac{1}{|B|^{\lambda_{1}}\|\chi_{B}\|_{\vec{p_{1}}}}
\sup\limits_{y\in B}|b(y)-b_{B}|\|\chi_{B}\|_{\vec{p_{1}}}
\end{align*}
\begin{align*}
&\leq\frac{C}{|B|^{\lambda_{1}}}\frac{\|(b-b_{B})\chi_{B}\|_{L^{1}(\mathbb{R}^{n})}}{\|\chi_{B}\|_{L^{1}(\mathbb{R}^{n})}}\\
&\leq\frac{C}{|B|^{\lambda_{1}}\|\chi_{B}\|_{\vec{p}}}\|(b-b_{B})\chi_{B}\|_{\vec{p}}.
\end{align*}
To deal with the above term, we note that
\begin{align*}
\|(b-b_{B})\chi_{B}\|_{\vec{p}}&\leq\Big\|\big(b(\cdot)-\frac{1}{|B|}\int_{B}b(z)dz\big)\chi_{B}(\cdot)\Big\|_{\vec{p}}\\
&=\Big\|\frac{\chi_{B}(\cdot)}{|B|}\int_{B}\big(b(\cdot)-b(z)\big)dz\Big\|_{\vec{p}}\\
&\leq\Big\|\frac{\chi_{B}(\cdot)}{|B|}\int_{|z|<|\cdot|}\big(b(\cdot)-b(z)\big)\chi_{B}(z)dz\Big\|_{\vec{p}}\\
&\quad+\Big\|\frac{\chi_{B}(\cdot)}{|B|}\int_{|z|\geq|\cdot|}\big(b(\cdot)-b(z)\big)\chi_{B}(z)dz\Big\|_{\vec{p}}\\
&=:G_{1}+G_{2}.
\end{align*}
The $(\mathcal{B}^{\vec{p_{2}},\lambda_{2}}(\mathbb{R}^n),\mathcal{B}^{\vec{p},\lambda}(\mathbb{R}^n))$ boundedness of $\mathcal{H}_{b}$ allows us to estimate $G_{1}$ as
\begin{align*}
G_{1}&\leq\Big\|\frac{\chi_{B}(\cdot)}{|B|}|\cdot|^{n}\frac{1}{|\cdot|^{n}}\int_{|z|<|\cdot|}\big(b(\cdot)-b(z)\big)\chi_{B}(z)dz\Big\|_{\vec{p}}\\
&\leq C\Big\|\chi_{B}(\cdot)\mathcal{H}_{b}(\chi_{B})(\cdot)\Big\|_{\vec{p}}\\
&\leq C|B|^{\lambda}\|\chi_{B}\|_{\vec{p}}\|\mathcal{H}_{b}(\chi_{B})\|_{\mathcal{B}^{\vec{p},\lambda}(\mathbb{R}^{n})}\\
&\leq C|B|^{\lambda}\|\chi_{B}\|_{\vec{p}}\|\chi_{B}\|_{\mathcal{B}^{\vec{p_{2}},\lambda_{2}}(\mathbb{R}^{n})}\\
&\leq C|B|^{\lambda_{1}}\|\chi_{B}\|_{\vec{p}}.
\end{align*}
By the $(\mathcal{B}^{\vec{p_{2}},\lambda_{2}}(\mathbb{R}^n),\mathcal{B}^{\vec{p},\lambda}(\mathbb{R}^n))$ boundedness of $\mathcal{H}^{*}_{b}$, it is easy to check that
\begin{align*}
G_{2}&\leq\Big\|\frac{\chi_{B}(\cdot)}{|B|}\int_{|z|\geq|\cdot|}\frac{\big(b(\cdot)-b(z)\big)\chi_{B}(z)}{|z|^{n}}|z|^{n}dz\Big\|_{\vec{p}}\\
&\leq C\Big\|\chi_{B}(\cdot)\mathcal{H}^{*}_{b}(\chi_{B})(\cdot)\Big\|_{\vec{p}}\\
&\leq C|B|^{\lambda}\|\chi_{B}\|_{\vec{p}}\|\mathcal{H}^{*}_{b}(\chi_{B})\|_{\mathcal{B}^{\vec{p},\lambda}(\mathbb{R}^{n})}\\
&\leq C|B|^{\lambda}\|\chi_{B}\|_{\vec{p}}\|\chi_{B}\|_{\mathcal{B}^{\vec{p_{2}},\lambda_{2}}(\mathbb{R}^{n})}\\
&\leq C|B|^{\lambda_{1}}\|\chi_{B}\|_{\vec{p}},
\end{align*}
where we have used the condition
$$\lambda=\lambda_{1}+\lambda_{2}.$$

We thus have established the following inequality if we combine the above estimates for $G_{1}$ and $G_{2}$,
$$\frac{\|(b-b_{B})\chi_{B}\|_{\vec{p_{1}}}}{|B|^{\lambda_{1}}\|\chi_{B}\|_{\vec{p_{1}}}}\leq\frac{C}{|B|^{\lambda_{1}}\|\chi_{B}\|_{\vec{p}}}
(|B|^{\lambda_{1}}\|\chi_{B}\|_{\vec{p}})\leq C.$$
The proof of Theorem 1.4 is finished.
\end{proof}

Let us now give the proof of Theorem 1.5.

\begin{proof}[\bf Proof of Theorem 1.5]

The proof of Theorem 1.5 is based on the following key lemma. For this purpose, we start with this necessary lemma.

\begin{lemma}
Let $1<\vec{p}<\infty$, $-\frac{1}{n}\sum^{n}_{i=1}\frac{1}{p_{i}}<\lambda<0$ and $i, k\in\mathbb{Z}$. If $b\in\mathfrak{C}^{\vec{p},\lambda}(\mathbb{R}^{n})$, then
$$|b(y)-b_{B_{k}}|\leq|b(y)-b_{B_{j}}|+C\max{\{|B_{k}|^{\lambda},|B_{j}|^{\lambda}\}}\|b\|_{\mathfrak{C}^{\vec{p},\lambda}(\mathbb{R}^{n})}.$$
\end{lemma}
${\bf Proof.}$
Using the H\"{o}lder inequality on mixed Lebesgue space to $\vec{p}$ and $\vec{p}'$, we have
\begin{align*}
|b_{B_{i}}-b_{B_{i+1}}|&\leq\frac{1}{|B_{i}|}\int_{B_{i}}|b(y)-b_{B_{i+1}}|dy\\
&\leq\frac{1}{|B_{i}|}\|\chi_{B_{i+1}}(b-b_{B_{i+1}})\|_{1}\\
&\leq\frac{1}{|B_{i}|}\|\chi_{B_{i+1}}\|_{\vec{p}'}\|(b-b_{B_{i+1}})\chi_{B_{i+1}}\|_{\vec{p}}\\
&\leq\frac{|B_{i+1}|^{\lambda}}{|B_{i}|}\|\chi_{B_{i+1}}\|_{\vec{p}'}\|\chi_{B_{i+1}}\|_{\vec{p}}\|b\|_{\mathfrak{C}^{\vec{p},\lambda}(\mathbb{R}^{n})}\\
&\leq C|B_{i+1}|^{\lambda}\|b\|_{\mathfrak{C}^{\vec{p},\lambda}(\mathbb{R}^{n})}.
\end{align*}

If $k<j$, then
\begin{align*}
|b(y)-b_{B_{k}}|&\leq|b(y)-b_{B_{j}}|+|b_{B_{k}}-b_{B_{j}}|\\
&\leq|b(y)-b_{B_{j}}|+\sum^{j-1}_{i=k}|b_{B_{i}}-b_{B_{i+1}}|\\
&\leq|b(y)-b_{B_{j}}|+C\sum^{j-1}_{i=k}|B_{i+1}|^{\lambda}\|b\|_{\mathfrak{C}^{\vec{p},\lambda}(\mathbb{R}^{n})}\\
&\leq|b(y)-b_{B_{j}}|+C|B_{k}|^{\lambda}\|b\|_{\mathfrak{C}^{\vec{p},\lambda}(\mathbb{R}^{n})}.
\end{align*}
If $j<k$, then

\begin{align*}
|b(y)-b_{B_{k}}|&\leq|b(y)-b_{B_{j}}|+\sum^{k-1}_{i=j}|b_{B_{i}}-b_{B_{i+1}}|\\
&\leq|b(y)-b_{B_{j}}|+C\sum^{k-1}_{i=j}|B_{i+1}|^{\lambda}\|b\|_{\mathfrak{C}^{\vec{p},\lambda}(\mathbb{R}^{n})}\\
&\leq|b(y)-b_{B_{j}}|+C|B_{j}|^{\lambda}\|b\|_{\mathfrak{C}^{\vec{p},\lambda}(\mathbb{R}^{n})}.
\end{align*}
The proof of Lemma 2.1 is completed.

Having disposed of this preliminary step, we can now return to the proof of Theorem 1.5.

$(a)\Rightarrow(b)$ Let $f$ be a function in $\mathcal{B}^{\vec{p},\lambda}(\mathbb{R}^n)$, $b\in\mathfrak{C}^{\vec{p},\lambda}(\mathbb{R}^{n})$. For a fixed ball $B(0,r)=B_{k_{0}}$ with $k_{0}\in\mathbb{Z}$, the following inequalities are true
\begin{equation} \label{2.3}
\|(\mathcal{H}_{b}f)\chi_{B_{k_{0}}}\|_{\vec{p}}\lesssim|B_{k_{0}}|^{2\lambda}
\|\chi_{B_{k_{0}}}\|_{\vec{p}}\|b\|_{\mathfrak{C}^{\vec{p},\lambda}(\mathbb{R}^{n})}
\|f\|_{\mathcal{B}^{\vec{p},\lambda}(\mathbb{R}^{n})}
\end{equation}
and
\begin{equation} \label{2.4}
\|(\mathcal{H}^{*}_{b}f)\chi_{B_{k_{0}}}\|_{\vec{p}}\lesssim|B_{k_{0}}|^{2\lambda}\|\chi_{B_{k_{0}}}\|_{\vec{p}}
\|b\|_{\mathfrak{C}^{\vec{p},\lambda}(\mathbb{R}^{n})}
\|f\|_{\mathcal{B}^{\vec{p},\lambda}(\mathbb{R}^{n})}.
\end{equation}

Let's dispose of (\ref{2.3}) now. We write
\begin{align*}
\|\mathcal{H}_{b}f(\cdot)\chi_{B_{k_{0}}}(\cdot)\|_{\vec{p}}&=
\Big\|\chi_{B_{k_{0}}}(\cdot)\frac{1}{|\cdot|^{n}}\int_{|y|<|\cdot|}\big(b(\cdot)-b(y)\big)f(y)dy\Big\|_{\vec{p}}\\
&\leq\sum^{k_{0}}_{k=-\infty}\Big\|\chi_{C_{k}}(\cdot)\frac{1}{|\cdot|^{n}}\int_{B_{k}}\big(b(\cdot)-b(y)\big)f(y)dy\Big\|_{\vec{p}}\\
&\lesssim\sum^{k_{0}}_{k=-\infty}\Big\|\chi_{C_{k}}(\cdot)\frac{1}{|\cdot|^{n}}\sum^{k}_{i=-\infty}
\int_{C_{i}}\big(b(\cdot)-b_{B_{k}}\big)f(y)dy\Big\|_{\vec{p}}\\
&\quad+\sum^{k_{0}}_{k=-\infty}\Big\|\chi_{C_{k}}(\cdot)\frac{1}{|\cdot|^{n}}\sum^{k}_{i=-\infty}
\int_{C_{i}}\big(b(y)-b_{B_{k}}\big)f(y)dy\Big\|_{\vec{p}}\\
&=:J+JJ.
\end{align*}
Here we give the calculation of $J$. For $1=\frac{1}{\vec{p}}+\frac{1}{\vec{p}'}$, by the H\"{o}lder inequality on mixed Lebesgue space, we get
\begin{align*}
J&\lesssim\sum^{k_{0}}_{k=-\infty}\frac{1}{|B_{k}|}\Big\|\chi_{C_{k}}\big(b-b_{B_{k}}\big)\sum^{k}_{i=-\infty}\int_{C_{i}}f(y)dy\Big\|_{\vec{p}}\\
&\leq\sum^{k_{0}}_{k=-\infty}\frac{1}{|B_{k}|}\Big\|\chi_{B_{k}}\big(b-b_{B_{k}}\big)\Big\|_{\vec{p}}
\sum^{k}_{i=-\infty}\int_{B_{i}}|f(y)|dy\\
&\leq\sum^{k_{0}}_{k=-\infty}\frac{1}{|B_{k}|}|B_{k}|^{\lambda}\big\|\chi_{B_{k}}\big\|_{\vec{p}}\big\|b\big\|_{\mathfrak{C}^{\vec{p},\lambda}(\mathbb{R}^{n})}
\sum^{k}_{i=-\infty}\big\|f\chi_{B_{i}}\big\|_{\vec{p}}\big\|\chi_{B_{i}}\big\|_{\vec{p}'}\\
&\leq\|b\|_{\mathfrak{C}^{\vec{p},\lambda}(\mathbb{R}^{n})}\|f\|_{\mathcal{B}^{\vec{p},\lambda}(\mathbb{R}^{n})}
\sum^{k_{0}}_{k=-\infty}|B_{k}|^{\lambda-1}\big\|\chi_{B_{k}}\big\|_{\vec{p}}\sum^{k}_{i=-\infty}|B_{i}|^{\lambda+1}\\
&\lesssim|B_{k_{0}}|^{2\lambda}\|\chi_{B_{k_{0}}}\|_{\vec{p}}
\|b\|_{\mathfrak{C}^{\vec{p},\lambda}(\mathbb{R}^{n})}\|f\|_{\mathcal{B}^{\vec{p},\lambda}(\mathbb{R}^{n})},
\end{align*}
we have used the conditions $-\frac{1}{n}\sum^{n}_{i=1}\frac{1}{2p_{i}}<\lambda<0$.

Using Lemma 2.1, we have
\begin{align*}
JJ&\lesssim\sum^{k_{0}}_{k=-\infty}\frac{1}{|B_{k}|}\Big\|\chi_{C_{k}}\sum^{k}_{i=-\infty}\int_{C_{i}}
\big|b(y)-b_{B_{k}}\big|\big|f(y)\big|dy\Big\|_{\vec{p}}\\
&\lesssim\sum^{k_{0}}_{k=-\infty}\frac{1}{|B_{k}|}\Big\|\chi_{B_{k}}\sum^{k}_{i=-\infty}\int_{B_{i}}
\big|b(y)-b_{B_{i}}\big|\big|f(y)\big|dy\Big\|_{\vec{p}}\\
&\quad+C\|b\|_{\mathfrak{C}^{\vec{p},\lambda}(\mathbb{R}^{n})}\sum^{k_{0}}_{k=-\infty}\frac{1}{|B_{k}|}\Big\|\chi_{B_{k}}\sum^{k}_{i=-\infty}\int_{B_{i}}
|B_{i}|^{\lambda}|f(y)|dy\Big\|_{\vec{p}}\\
&=:JJ_{1}+JJ_{2}.
\end{align*}
We may choose $\vec{s}$ such that $\frac{1}{\vec{p}'}=\frac{1}{\vec{p}}+\frac{1}{\vec{s}}$, then $\sum^{n}_{i=1}\frac{1}{s_{i}}=\sum^{n}_{i=1}\frac{1}{p_{i}'}-\sum^{n}_{i=1}\frac{1}{p_{i}}$. Repeated application of H\"{o}lder's inequality on mixed Lebesgue space shows that
\begin{align*}
JJ_{1}&\leq\sum^{k_{0}}_{k=-\infty}\frac{1}{|B_{k}|}\big\|\chi_{B_{k}}\big\|_{\vec{p}}
\sum^{k}_{i=-\infty}\big\|(b-b_{B_{i}})\chi_{B_{i}}\big\|_{\vec{p}'}\big\|f\chi_{B_{i}}\big\|_{\vec{p}}\\
&\leq\|f\|_{\mathcal{B}^{\vec{p},\lambda}(\mathbb{R}^{n})}\sum^{k_{0}}_{k=-\infty}\frac{1}{|B_{k}|}\big\|\chi_{B_{k}}\big\|_{\vec{p}}
\sum^{k}_{i=-\infty}\big\|(b-b_{B_{i}})\chi_{B_{i}}\big\|_{\vec{p}}\big\|\chi_{B_{i}}\big\|_{\vec{s}}|B_{i}|^{\lambda}\big\|\chi_{B_{i}}\big\|_{\vec{p}}\\
&\leq\|b\|_{\mathfrak{C}^{\vec{p},\lambda}(\mathbb{R}^{n})}\|f\|_{\mathcal{B}^{\vec{p},\lambda}(\mathbb{R}^{n})}
\sum^{k_{0}}_{k=-\infty}\frac{1}{|B_{k}|}\big\|\chi_{B_{k}}\big\|_{\vec{p}}
\sum^{k}_{i=-\infty}\big\|\chi_{B_{i}}\big\|_{\vec{p}}\big\|\chi_{B_{i}}\big\|_{\vec{p}'}|B_{i}|^{2\lambda}\\
&\leq\|b\|_{\mathfrak{C}^{\vec{p},\lambda}(\mathbb{R}^{n})}\|f\|_{\mathcal{B}^{\vec{p},\lambda}(\mathbb{R}^{n})}
\sum^{k_{0}}_{k=-\infty}\frac{1}{|B_{k}|}\big\|\chi_{B_{k}}\big\|_{\vec{p}}\sum^{k}_{i=-\infty}|B_{i}|^{2\lambda+1}\\
&\leq|B_{k_{0}}|^{2\lambda}\|\chi_{B_{k_{0}}}\|_{\vec{p}}
\|b\|_{\mathfrak{C}^{\vec{p},\lambda}(\mathbb{R}^{n})}\|f\|_{\mathcal{B}^{\vec{p},\lambda}(\mathbb{R}^{n})},
\end{align*}
and the fact $-\frac{1}{n}\sum^{n}_{i=1}\frac{1}{2p_{i}}<\lambda<0$ has been used.

For $-\frac{1}{n}\sum^{n}_{i=1}\frac{1}{2p_{i}}<\lambda<0$, by the H\"{o}lder inequality on mixed Lebesgue space, the term $JJ_{2}$ can be bounded by
\begin{align*}
JJ_{2}&\leq\|b\|_{\mathfrak{C}^{\vec{p},\lambda}(\mathbb{R}^{n})}\sum^{k_{0}}_{k=-\infty}\frac{1}{|B_{k}|}\Big\|\chi_{B_{k}}\Big\|_{\vec{p}}
\sum^{k}_{i=-\infty}|B_{i}|^{\lambda}\int_{B_{i}}|f(y)|dy\\
&\leq\|b\|_{\mathfrak{C}^{\vec{p},\lambda}(\mathbb{R}^{n})}\sum^{k_{0}}_{k=-\infty}\frac{1}{|B_{k}|}\Big\|\chi_{B_{k}}\Big\|_{\vec{p}}
\sum^{k}_{i=-\infty}|B_{i}|^{\lambda}\big\|f\chi_{B_{i}}\big\|_{\vec{p}}\big\|\chi_{B_{i}}\big\|_{\vec{p}'}\\
&\leq\|b\|_{\mathfrak{C}^{\vec{p},\lambda}(\mathbb{R}^{n})}\|f\|_{\mathcal{B}^{\vec{p},\lambda}(\mathbb{R}^{n})}
\sum^{k_{0}}_{k=-\infty}\frac{1}{|B_{k}|}\big\|\chi_{B_{k}}\big\|_{\vec{p}}\sum^{k}_{i=-\infty}|B_{i}|^{2\lambda+1}\\
&\leq|B_{k_{0}}|^{2\lambda}\|\chi_{B_{k_{0}}}\|_{\vec{p}}
\|b\|_{\mathfrak{C}^{\vec{p},\lambda}(\mathbb{R}^{n})}\|f\|_{\mathcal{B}^{\vec{p},\lambda}(\mathbb{R}^{n})}.
\end{align*}
We have thus proved the estimate for $JJ$. (\ref{2.3}) is based on the above estimates for $J$ and $JJ$.

With a slight modification of the proofs for (\ref{2.2}) and (\ref{2.3}), (\ref{2.4}) can be obtained easily, we omit its proof here for the similarity.

$(b)\Rightarrow(a)$
In this case, suppose that $\mathcal{H}_{b}$ and $\mathcal{H}^{*}_{b}$ are bounded operators from $\mathcal{B}^{\vec{p},\lambda}(\mathbb{R}^n)$ to $\mathcal{B}^{\vec{p},2\lambda}(\mathbb{R}^n)$, we will show that
$$b\in\mathfrak{C}^{\vec{p},\lambda}(\mathbb{R}^{n}).$$
We only need to show that for a fixed ball $B=B(0,r)$,
$$\frac{\|(b-b_{B})\chi_{B}\|_{\vec{p}}}{|B|^{\lambda}\|\chi_{B}\|_{\vec{p}}}\leq C.$$
To deal with the above term, it will be necessary to note that
\begin{align*}
\frac{\|(b-b_{B})\chi_{B}\|_{\vec{p}}}{|B|^{\lambda}\|\chi_{B}\|_{\vec{p}}}
&\leq\frac{1}{|B|^{\lambda}\|\chi_{B}\|_{\vec{p}}}\Big\|\big(b(\cdot)-\frac{1}{|B|}\int_{B}b(z)dz\big)\chi_{B}(\cdot)\Big\|_{\vec{p}}\\
&=\frac{1}{|B|^{\lambda}\|\chi_{B}\|_{\vec{p}}}\Big\|\frac{\chi_{B}(\cdot)}{|B|}\int_{B}\big(b(\cdot)-b(z)\big)dz\Big\|_{\vec{p}}
\end{align*}
\begin{align*}
&\leq\frac{1}{|B|^{\lambda}\|\chi_{B}\|_{\vec{p}}}\Big\|\frac{\chi_{B}(\cdot)}{|B|}\int_{|z|<|\cdot|}\big(b(\cdot)-b(z)\big)\chi_{B}(z)dz\Big\|_{\vec{p}}\\
&\quad+\frac{1}{|B|^{\lambda}\|\chi_{B}\|_{\vec{p}}}\Big\|\frac{\chi_{B}(\cdot)}{|B|}\int_{|z|\geq|\cdot|}\big(b(\cdot)-b(z)\big)\chi_{B}(z)dz\Big\|_{\vec{p}}\\
&=:K_{1}+K_{2}.
\end{align*}
The $(\mathcal{B}^{\vec{p},\lambda}(\mathbb{R}^n),\mathcal{B}^{\vec{p},2\lambda}(\mathbb{R}^n))$ boundedness of $\mathcal{H}_{b}$ produces the following estimate for the term $K_{1}$,
\begin{align*}
K_{1}&\leq\frac{1}{|B|^{\lambda}\|\chi_{B}\|_{\vec{p}}}\Big\|\frac{\chi_{B}(\cdot)}{|B|}|\cdot|^{n}\frac{1}{|\cdot|^{n}}\int_{|z|<|\cdot|}\big(b(\cdot)-b(z)\big)\chi_{B}(z)dz\Big\|_{\vec{p}}\\
&\leq\frac{C}{|B|^{\lambda}\|\chi_{B}\|_{\vec{p}}}\Big\|\chi_{B}(\cdot)\mathcal{H}_{b}(\chi_{B})(\cdot)\Big\|_{\vec{p}}\\
&\leq C|B|^{\lambda}\|\mathcal{H}_{b}(\chi_{B})\|_{\mathcal{B}^{\vec{p},2\lambda}(\mathbb{R}^{n})}\\
&\leq C|B|^{\lambda}\|\chi_{B}\|_{\mathcal{B}^{\vec{p},\lambda}(\mathbb{R}^{n})}\\
&\leq C.
\end{align*}
By the $(\mathcal{B}^{\vec{p},\lambda}(\mathbb{R}^n),\mathcal{B}^{\vec{p},2\lambda}(\mathbb{R}^n))$ boundedness of $\mathcal{H}^{*}_{b}$, the following can be confirmed easily
\begin{align*}
K_{2}&\leq\frac{1}{|B|^{\lambda}\|\chi_{B}\|_{\vec{p}}}\Big\|\frac{\chi_{B}(\cdot)}{|B|}\int_{|z|\geq|\cdot|}\frac{\big(b(\cdot)-b(z)\big)\chi_{B}(z)}{|z|^{n}}|z|^{n}dz\Big\|_{\vec{p}}\\
&\leq \frac{C}{|B|^{\lambda}\|\chi_{B}\|_{\vec{p}}}\Big\|\chi_{B}(\cdot)\mathcal{H}^{*}_{b}(\chi_{B})(\cdot)\Big\|_{\vec{p}}\\
&\leq C|B|^{\lambda}\|\mathcal{H}^{*}_{b}(\chi_{B})\|_{\mathcal{B}^{\vec{p},2\lambda}(\mathbb{R}^{n})}\\
&\leq C|B|^{\lambda}\|\chi_{B}\|_{\mathcal{B}^{\vec{p},\lambda}(\mathbb{R}^{n})}\\
&\leq C.
\end{align*}
Thus, $K_{1}$ together with $K_{2}$ implies that

$$\frac{\|(b-b_{B})\chi_{B}\|_{\vec{p}}}{|B|^{\lambda}\|\chi_{B}\|_{\vec{p}}}\leq C.$$
This completes the proof of Theorem 1.5.
\end{proof}

\section{Behavior of  Hardy type operators on mixed $\lambda$-central Morrey spaces}

In this section, we will give some characterizations of the mixed central bounded mean oscillation spaces, and then establish some boundedness of Hardy type operators on mixed $\lambda$-central Morrey spaces. We formulate  main results as follows.

\begin{theorem}
Let $1<\vec{p}<\infty$, $-\frac{1}{n}\sum^{n}_{i=1}\frac{1}{p_{i}}<\lambda<0$, and $1=\frac{1}{\vec{p}}+\frac{1}{\vec{p}'}$. Then

(a)\quad $\mathcal{H}$ is a bounded operator from $\mathcal{B}^{\vec{p},\lambda}(\mathbb{R}^n)$ to $\mathcal{B}^{\vec{p},\lambda}(\mathbb{R}^n)$;

(b)\quad $\mathcal{H}^{*}$ is a bounded operator from $\mathcal{B}^{\vec{p}',\lambda}(\mathbb{R}^n)$ to $\mathcal{B}^{\vec{p}',\lambda}(\mathbb{R}^n)$.
\end{theorem}

\begin{theorem}
Let $\vec{p}, \lambda$, and $\vec{p}'$ be as in Theorem 3.1. Then the following statements are equivalent:

(a)\quad $b\in \mathrm{CMO}^{\vec{p}}(\mathbb{R}^{n})\cap \mathrm{CMO}^{\vec{p}'}(\mathbb{R}^{n})$;

(b)\quad $\mathcal{H}_{b}$ and $\mathcal{H}^{*}_{b}$ are bounded operators on both $\mathcal{B}^{\vec{p},\lambda}(\mathbb{R}^n)$ and $\mathcal{B}^{\vec{p}',\lambda}(\mathbb{R}^n)$.
\end{theorem}

Theorem 3.1 and Theorem 3.2 state that the boundedness of $\mathcal{H}$ and $\mathcal{H}_{b}$ on $\mathcal{B}^{\vec{p},\lambda}(\mathbb{R}^n)$ does not involve any further assumptions about the indexes $\lambda$ and $\vec{p}$ comparing to that of Theorem 1.4 and Theorem 1.5.

\hspace{-22pt}{\bf Remark 3.1} {The corresponding boundedness for $\mathcal{H}_{\alpha}$ and $\mathcal{H}_{\alpha,b}$ on $\mathcal{B}^{\vec{p}',\lambda}(\mathbb{R}^n)$ is also true.
}

The proofs of Theorem 3.1 and Theorem 3.2 based on the following lemmas about the boundedness of Hardy type operators on mixed Lebesgue spaces.

\begin{proof}[\bf Proof of Theorem 3.1] Following the notations of Section 2, we need to show that for a fixed ball $B(0.r)=B_{k_{0}}$ with $k_{0}\in\mathbb{Z}$, there exist constants $C>0$ such that
\begin{equation} \label{3.1}
\|(\mathcal{H}f)\chi_{B_{k_{0}}}\|_{\vec{p}}\leq C|B_{k_{0}}|^{\lambda}
\|\chi_{B_{k_{0}}}\|_{\vec{p}}\|f\|_{\mathcal{B}^{\vec{p},\lambda}(\mathbb{R}^{n})}
\end{equation}
and
\begin{equation} \label{3.2}
\|(\mathcal{H}^{*}f)\chi_{B_{k_{0}}}\|_{\vec{p}}\leq C|B_{k_{0}}|^{\lambda}\|\chi_{B_{k_{0}}}\|_{\vec{p}'}
\|f\|_{\mathcal{B}^{\vec{p}',\lambda}(\mathbb{R}^{n})}.
\end{equation}
For any ball $B$, let $2B=B(0,2r)$, we write $f$ as $f=f\chi_{2B_{k_{0}}}+f\chi_{(2B_{k_{0}})^{c}}=:f_{1}+f_{2}$, then
\begin{align*}
\|(\mathcal{H}f)\chi_{B_{k_{0}}}\|_{\vec{p}}&\leq\|(\mathcal{H}f_{1})\chi_{B_{k_{0}}}\|_{\vec{p}}+\|(\mathcal{H}f_{2})\chi_{B_{k_{0}}}\|_{\vec{p}}\\
&=:T_{1}+T_{2}.
\end{align*}

We first estimate $T_{1}$. By the boundedness of Hardy operator $\mathcal{H}$ on mixed Lebesgue space $L^{\vec{p}}(\mathbb{R}^{n})$ (see \cite{WM2}), we get
$$T_{1}\leq C\|f\chi_{2B_{k_{0}}}\|_{\vec{p}}\leq C|B_{k_{0}}|^{\lambda}
\|\chi_{B_{k_{0}}}\|_{\vec{p}}\|f\|_{\mathcal{B}^{\vec{p},\lambda}(\mathbb{R}^{n})}.$$

Now we turn to prove $T_{2}$. Using H\"{o}lder's inequality on mixed Lebesgue spaces, since
\begin{align*}
|\mathcal{H}f_{2}(x)|&\leq\frac{1}{|x|^{n}}\int_{|y|<|x|}|f_{2}(y)|dy\\
&\leq C\sum^{\infty}_{k=2k_{0}}\frac{1}{|B_{k}|}\int_{C_{k}}|f(y)|dy\\
&\leq C\sum^{\infty}_{k=2k_{0}}\frac{1}{|B_{k}|}\big\|f\chi_{B_{k}}\big\|_{\vec{p}}\big\|\chi_{B_{k}}\big\|_{\vec{p}'}\\
&\leq C\|f\|_{\mathcal{B}^{\vec{p},\lambda}(\mathbb{R}^{n})}\sum^{\infty}_{k=2k_{0}}|B_{k}|^{\lambda}\\
&\leq C\|f\|_{\mathcal{B}^{\vec{p},\lambda}(\mathbb{R}^{n})}|B_{k_{0}}|^{\lambda},
\end{align*}
where we have used the fact that $-\frac{1}{n}\sum^{n}_{i=1}\frac{1}{p_{i}}<\lambda<0.$
Therefore,
$$T_{2}\leq C|B_{k_{0}}|^{\lambda}
\|\chi_{B_{k_{0}}}\|_{\vec{p}}\|f\|_{\mathcal{B}^{\vec{p},\lambda}(\mathbb{R}^{n})}.$$

Combining the estimates of $T_{1}$ and $T_{2}$, (\ref{3.1}) is obtained. (\ref{3.2}) follows by the same method as that of (\ref{3.1}).
\end{proof}

\begin{proof}[\bf Proof of Theorem 3.2]

$(a)\Rightarrow(b)$ In this case, the task is to show that for a fixed ball $B(0,r)=B_{k_{0}}$ with $k_{0}\in\mathbb{Z}$, let $f$ be a function in $\mathcal{B}^{\vec{p},\lambda}(\mathbb{R}^n)$, $b\in\mathrm{CMO}^{\vec{p}}(\mathbb{R}^{n})\cap \mathrm{CMO}^{\vec{p}'}(\mathbb{R}^{n})$, there exist constants $C>0$ such that the following inequalities are true:
\begin{equation} \label{3.3}
\|(\mathcal{H}_{b}f)\chi_{B_{k_{0}}}\|_{\vec{p}}\lesssim|B_{k_{0}}|^{\lambda}
\|\chi_{B_{k_{0}}}\|_{\vec{p}}\|f\|_{\mathcal{B}^{\vec{p},\lambda}(\mathbb{R}^{n})},
\end{equation}

\begin{equation} \label{3.4}
\|(\mathcal{H}^{*}_{b}f)\chi_{B_{k_{0}}}\|_{\vec{p}}\lesssim|B_{k_{0}}|^{\lambda}\|\chi_{B_{k_{0}}}\|_{\vec{p}}
\|f\|_{\mathcal{B}^{\vec{p},\lambda}(\mathbb{R}^{n})},
\end{equation}

\begin{equation} \label{3.5}
\|(\mathcal{H}_{b}f)\chi_{B_{k_{0}}}\|_{\vec{p}'}\lesssim|B_{k_{0}}|^{\lambda}
\|\chi_{B_{k_{0}}}\|_{\vec{p}'}\|f\|_{\mathcal{B}^{\vec{p}',\lambda}(\mathbb{R}^{n})}
\end{equation}
and
\begin{equation} \label{3.6}
\|(\mathcal{H}^{*}_{b}f)\chi_{B_{k_{0}}}\|_{\vec{p}'}\lesssim|B_{k_{0}}|^{\lambda}\|\chi_{B_{k_{0}}}\|_{\vec{p}'}
\|f\|_{\mathcal{B}^{\vec{p}',\lambda}(\mathbb{R}^{n})}.
\end{equation}
Let's dispose of (\ref{3.3}) now. Decomposing $f=f\chi_{2B_{k_{0}}}+f\chi_{(2B_{k_{0}})^{c}}=:f_{1}+f_{2}$ derives that
\begin{align*}
\|(\mathcal{H}_{b}f)\chi_{B_{k_{0}}}\|_{\vec{p}}&\leq\|(\mathcal{H}_{b}f_{1})\chi_{B_{k_{0}}}\|_{\vec{p}}+\|(\mathcal{H}_{b}f_{2})\chi_{B_{k_{0}}}\|_{\vec{p}}\\
&=:TT_{1}+TT_{2}.
\end{align*}

We can infer from Theorem 1.3 that
\begin{align*}
TT_{1}&\leq\|(\mathcal{H}_{b}f_{1})\chi_{B_{k_{0}}}\|_{\vec{p}}\leq C\|f\chi_{2B_{k_{0}}}\|_{\vec{p}}\\
&\leq C|B_{k_{0}}|^{\lambda}\|\chi_{B_{k_{0}}}\|_{\vec{p}}\|f\|_{\mathcal{B}^{\vec{p},\lambda}(\mathbb{R}^{n})}.
\end{align*}
Next we estimate $TT_{2}$. By the definitions of $\mathcal{H}_{b}f$ and $\mathrm{CMO}^{\vec{p}}(\mathbb{R}^{n})$, $TT_{2}$ can be shown
\begin{align*}
TT_{2}&\leq\Big\|\chi_{B_{k_{0}}}(\cdot)\frac{1}{|\cdot|^{n}}\int_{|y|<|\cdot|}|b(\cdot)-b(y)||f_{2}(y)|dy\Big\|_{\vec{p}}\\
&\leq\Big\|\chi_{B_{k_{0}}}(\cdot)\sum^{\infty}_{k=2k_{0}}\frac{1}{|B_{k}|}\int_{B_{k}}|b(\cdot)-c||f(y)|dy\Big\|_{\vec{p}}\\
&\quad+\Big\|\chi_{B_{k_{0}}}(\cdot)\sum^{\infty}_{k=2k_{0}}\frac{1}{|B_{k}|}\int_{B_{k}}|b(y)-c||f(y)|dy\Big\|_{\vec{p}}\\
&=:TT_{21}+TT_{22}.
\end{align*}
Applying H\"{o}lder's inequality on mixed Lebesgue space to $\vec{p}$ and $\vec{p}'$ can get the estimate for $TT_{21}$ as
\begin{align*}
TT_{21}&\leq\big\|(b-c)\chi_{B_{k_{0}}}\big\|_{\vec{p}}\sum^{\infty}_{k=2k_{0}}\frac{1}{|B_{k}|}\int_{B_{k}}|f(y)|dy\\
&\leq C\|\chi_{B_{k_{0}}}\|_{\vec{p}}\|b\|_{\mathrm{CMO}^{\vec{p}}(\mathbb{R}^{n})}\sum^{\infty}_{k=2k_{0}}\frac{1}{|B_{k}|}
\|f\chi_{B_{k}}\|_{\vec{p}}\|\chi_{B_{k}}\|_{\vec{p}'}\\
&\leq C\|\chi_{B_{k_{0}}}\|_{\vec{p}}\|b\|_{\mathrm{CMO}^{\vec{p}}(\mathbb{R}^{n})}\|f\|_{\mathcal{B}^{\vec{p},\lambda}(\mathbb{R}^{n})}
\sum^{\infty}_{k=2k_{0}}|B_{k}|^{\lambda}\\
&\leq C|B_{k_{0}}|^{\lambda}\|\chi_{B_{k_{0}}}\|_{\vec{p}}\|b\|_{\mathrm{CMO}^{\vec{p}}(\mathbb{R}^{n})}\|f\|_{\mathcal{B}^{\vec{p},\lambda}(\mathbb{R}^{n})}.
\end{align*}
The same result can be acquired for the term $TT_{22}$ as
\begin{align*}
TT_{22}&\leq C\big\|\chi_{B_{k_{0}}}\big\|_{\vec{p}}\sum^{\infty}_{k=2k_{0}}\frac{1}{|B_{k}|}\int_{B_{k}}|b(y)-c||f(y)|dy\\
&\leq C\|\chi_{B_{k_{0}}}\|_{\vec{p}}\sum^{\infty}_{k=2k_{0}}\frac{1}{|B_{k}|}\|(b-c)\chi_{B_{k}}\|_{\vec{p}'}\|f\chi_{B_{k}}\|_{\vec{p}}\\
&\leq C\|\chi_{B_{k_{0}}}\|_{\vec{p}}\|b\|_{\mathrm{CMO}^{\vec{p}'}(\mathbb{R}^{n})}\|f\|_{\mathcal{B}^{\vec{p},\lambda}(\mathbb{R}^{n})}
\sum^{\infty}_{k=2k_{0}}|B_{k}|^{\lambda}\\
&\leq C|B_{k_{0}}|^{\lambda}\|\chi_{B_{k_{0}}}\|_{\vec{p}}\|b\|_{\mathrm{CMO}^{\vec{p}'}(\mathbb{R}^{n})}\|f\|_{\mathcal{B}^{\vec{p},\lambda}(\mathbb{R}^{n})}.
\end{align*}

Combining the above estimates of $TT_{1}$ and $TT_{2}$, (\ref{3.3}) is shown. Similarly, we can prove (\ref{3.4}), (\ref{3.5}) and (\ref{3.6}).

$(b)\Rightarrow(a)$ Conversely, suppose that $\mathcal{H}_{b}$ and $\mathcal{H}^{*}_{b}$ are bounded operators on both $\mathcal{B}^{\vec{p},\lambda}(\mathbb{R}^n)$ and $\mathcal{B}^{\vec{p}',\lambda}(\mathbb{R}^n)$, we will show that
$$b\in \mathrm{CMO}^{\vec{p}}(\mathbb{R}^{n})\cap \mathrm{CMO}^{\vec{p}'}(\mathbb{R}^{n}).$$

Actually, we consider only the situation when $b\in \mathrm{CMO}^{\vec{p}}(\mathbb{R}^{n})$.  For the case of all $b\in \mathrm{CMO}^{\vec{p}'}(\mathbb{R}^{n})$, we know that both $\mathcal{H}_{b}$ and $\mathcal{H}^{*}_{b}$ map $L^{\vec{p}'}(\mathbb{R}^{n})$ into $L^{\vec{p}'}(\mathbb{R}^{n})$ via Remark 3.1 (see \cite{LZ3}). Therefore, a similar procedure can derive the desired result.

For any fixed $r>0$, denote $B(0,r)$ by $B$. Since $\mathcal{H}_{b}$ and $\mathcal{H}^{*}_{b}$ are bounded from $L^{\vec{p}}(\mathbb{R}^{n})$ to $L^{\vec{p}}(\mathbb{R}^{n})$, we have
\begin{align*}
\frac{\|(b-b_{B(0,r)})\chi_{B(0,r)}\|_{\vec{p}}}{\|\chi_{B(0,r)}\|_{\vec{p}}}&=
\frac{1}{\|\chi_{B}\|_{\vec{p}}}\Big\|\big(b(\cdot)-\frac{1}{|B|}\int_{B}b(z)dz\big)\chi_{B}(\cdot)\Big\|_{\vec{p}}\\
&=\frac{1}{\|\chi_{B}\|_{\vec{p}}}\Big\|\frac{\chi_{B}(\cdot)}{|B|}\int_{B}\big(b(\cdot)-b(z)\big)dz\Big\|_{\vec{p}}\\
&\leq\frac{1}{\|\chi_{B}\|_{\vec{p}}}\Big\|\frac{\chi_{B}(\cdot)}{|B|}\int_{|z|<|\cdot|}\big(b(\cdot)-b(z)\big)\chi_{B}(z)dz\Big\|_{\vec{p}}\\
&\quad+\frac{1}{\|\chi_{B}\|_{\vec{p}}}\Big\|\frac{\chi_{B}(\cdot)}{|B|}\int_{|z|\geq|\cdot|}\big(b(\cdot)-b(z)\big)\chi_{B}(z)dz\Big\|_{\vec{p}}\\
&\leq\frac{1}{\|\chi_{B}\|_{\vec{p}}}\Big\|\frac{\chi_{B}(\cdot)}{|B|}|\cdot|^{n}\frac{1}{|\cdot|^{n}}\int_{|z|<|\cdot|}\big(b(\cdot)-b(z)\big)\chi_{B}(z)dz\Big\|_{\vec{p}}\\
&\quad+\frac{1}{\|\chi_{B}\|_{\vec{p}}}\Big\|\frac{\chi_{B}(\cdot)}{|B|}\int_{|z|\geq|\cdot|}\frac{\big(b(\cdot)-b(z)\big)\chi_{B}(z)}{|z|^{n}}|z|^{n}dz\Big\|_{\vec{p}}\\
&\leq\frac{C}{\|\chi_{B}\|_{\vec{p}}}\Big\|\chi_{B}(\cdot)\mathcal{H}_{b}(\chi_{B})(\cdot)\Big\|_{\vec{p}}
+\frac{C}{\|\chi_{B}\|_{\vec{p}}}\Big\|\chi_{B}(\cdot)\mathcal{H}^{*}_{b}(\chi_{B})(\cdot)\Big\|_{\vec{p}}\\
&\leq\frac{C}{\|\chi_{B}\|_{\vec{p}}}\|\chi_{B}\|_{\vec{p}}\\
&\leq C,
\end{align*}
which finishes the proof of Theorem 3.2.
\end{proof}

\noindent$\textbf{Declarations.}$\\
The authors declare that there is no conflict of interests regarding the publication of this paper.
\noindent$\textbf{Acknowledgement.}$ \\
This research is supported by the National Natural Science Foundation of China (No. 12061069), The authors would like to express their gratitude to the referee for his/her very valuable comments.

\bigskip

\noindent Wenna Lu

\medskip

\noindent College of Mathematics and System Sciences,\ Xinjiang University\\
\"{U}r\"{u}mqi  830046,\ People's Republic of China

\smallskip

\noindent{\it E-mail address}: \texttt{luwnmath@126.com}

\bigskip

\noindent Jiang Zhou$^\ast$ (Corresponding author)

\medskip

\noindent College of Mathematics and System Sciences,\ Xinjiang University\\
\"{U}r\"{u}mqi  830046,\ People's Republic of China

\smallskip
\noindent{\it E-mail address}: \texttt{zhoujiang@xju.edu.cn}

\end{document}